\documentclass[12pt]{article}

\usepackage{mathrsfs}
\usepackage{color}
\usepackage{amssymb, amsthm, amsfonts, amsxtra, amsmath}
\usepackage{mathabx}

\usepackage{latexsym}
\usepackage[all]{xy}
\usepackage{graphics}
\usepackage{latexsym}
\usepackage{makeidx}
\usepackage{rotating}

\theoremstyle{change}

\newtheorem{Theorem}{Theorem}[section]
\newtheorem{Def}[Theorem]{Definition}
\newtheorem{Lem}[Theorem]{Lemma}
\newtheorem{Prop}[Theorem]{Proposition}

\newtheorem{Not}[Theorem]{Notation}

\date{}

\begin{document}

\hyphenation{Wo-ro-no-wicz}

\title{On a Morita equivalence between the duals of quantum $SU(2)$ and quantum $\widetilde{E}(2)$}
\author{Kenny De Commer\footnote{Supported in part by the ERC Advanced Grant 227458
OACFT ``Operator Algebras and Conformal Field Theory" }\\ \small Dipartimento di Matematica,  Universit\`{a} degli Studi di Roma Tor Vergata\\
\small Via della Ricerca Scientifica 1, 00133 Roma, Italy\\ \\ \small e-mail: decommer@mat.uniroma2.it}
\maketitle

\newcommand{\acnabla}{\nabla\!\!\!{^\shortmid}}
\newcommand{\undersetmin}[2]{{#1}\underset{\textrm{min}}{\otimes}{#2}}
\newcommand{\otimesud}[2]{\overset{#2}{\underset{#1}{\otimes}}}
\newcommand{\qbin}[2]{\left[ \begin{array}{c} #1 \\ #2 \end{array}\right]_{q^2}}

\newcommand{\otimesmin}{\underset{\textrm{min}}{\otimes}}
\newcommand{\bigback}{\!\!\!\!\!\!\!\!\!\!\!\!\!\!\!\!\!\!\!\!\!\!\!\!}

\abstract{\noindent Let $SU_q(2)$ and $\widetilde{E}_q(2)$ be Woronowicz' $q$-deformations of respectively the compact Lie group $SU(2)$ and the non-trivial double cover of the Lie group $E(2)$ of Euclidian transformations of the plane. We prove that, in some sense, their duals are `Morita equivalent locally compact quantum groups'. In more concrete terms, we prove that the von Neumann algebraic quantum groups $\mathscr{L}^{\infty}(SU_q(2))$ and $\mathscr{L}^{\infty}(\widetilde{E}_q(2))$ are unitary cocycle deformations of each other.}\vspace{0.3cm}

\noindent \emph{Keywords}: locally compact quantum group; von Neumann algebraic quantum group; quantization of classical Lie groups; cocycle twisting; Morita equivalence\\

\noindent AMS 2000 \emph{Mathematics subject classification}: 20G42; 46L65; 81R50; 16W30

%17B37: Quantum groups, quantized enveloping algebras
%20G42: quantized function algebras
%46L65: Functional analysis, deformations, quantizations
%81R50: Quantum groups and related algebraic methods
%16W30: Hopf algebras, co-algebras

\section*{Introduction}

\noindent This is part of a series of papers devoted to an intriguing correspondence between the quantizations of $SU(2)$, $\widetilde{E}(2)$ and $\widetilde{SU}(1,1)$, where the latter denotes the normalizer of $SU(1,1)$ inside $SL(2,\mathbb{C})$. In a sense, their duals form a trinity of `Morita equivalent locally compact quantum groups'. There then exists a `linking quantum groupoid' combining these three quantum groups into one global structure, and it is important to understand for example the (co)representation theory of this object.\\

\noindent In this paper, we will treat the `group von Neumann algebra of the linking quantum groupoid between the dual of $SU_q(2)$ and the dual of $\widetilde{E}_q(2)$'. We also treat part of its associated infinitesimal description. One will see that it bears some similarities to the `contraction procedure' which has been studied on the algebraic level in a series of papers by Celeghini and collaborators (see e.g. \cite{Cel1}), and in a C$^*$-algebraic framework by Woronowicz in \cite{Wor6}. But our philosophy is different: while those authors consider the passage between the two quantum groups as a kind of limit procedure, we construct a concrete object linking them. In a sense, we construct a bridge to cross the water, while in the other approach, one searches a place where the river is shallow enough to cross. In any case, we will try to comment at the appropriate places when there is a concrete resemblance between our theory and the contraction procedure. We also remark that a close connection between $SU_q(2)$, $\widetilde{E}_q(2)$ and $SU_q(1,1)$ is known in relation to the $q$-analogue of the Askey-Wilson function transform scheme (\cite{Koe2}, section 7).\\

\noindent The main tool in this paper is the theory of projective representations for locally compact quantum groups, developed in the final part of  \cite{DeC3} (and based on observations by A. Wassermann in \cite{Was1}). Indeed, using this theory, we showed in \cite{DeC1} (chapter 10) that whenever a compact quantum group has an infinite-dimensional irreducible projective representation, it allows for an `exotic deformation' into a non-compact locally compact quantum group. It turns out that there are (at least) two distinct such irreducible representations for $SU_q(2)$, which even have a quantum geometrical description: the first is associated to the homogeneous action of $SU_q(2)$ on the standard Podle\'{s} sphere, the other with its homogeneous action on the quantum projective plane (\cite{Haj1}). These will be related to the deformations of $SU_q(2)$ into respectively $\widetilde{E}_q(2)$ and $\widetilde{SU}_q(1,1)$.\\

\noindent The contents of this paper are as follows.\\

\noindent In \emph{the first section} we recall the analytic notions of von Neumann algebraic quantum groups (see \cite{Kus1},\cite{Kus2},\cite{VDae4}), von Neumann algebraic linking quantum groupoids (\cite{DeC1}) and unitary projective representations for (locally) compact quantum groups (\cite{DeC3},\cite{DeC1}), and the algebraic notions of bi-Galois objects (\cite{Sch1}) and co-linking weak Hopf algebras (\cite{Bic1}, \cite{DeC1}). We also recall the definitions of $SU_q(2)$ and $\widetilde{E}_q(2)$ on the Hopf $^*$-algebra level, and of the dual quantized universal enveloping algebras $U_q(su(2))$ and $U_q(e(2))$.\\

\noindent In \emph{the second section}, we begin by observing that the von Neumann algebraic completion of the action of $SU_q(2)$ on the standard Podle\'{s} sphere provides a projective unitary representation of $SU_q(2)$. We state the fact, proven in \cite{DeC4}, that the standard Podle\'{s} sphere can also be seen as a subquotient $^*$-algebra of a $^*$-Galois object for $U_q(su(2))$, with the associated infinitesimal action coming from the Miyashita-Ulbrich action on the Galois object. \\

\noindent In the \emph{third section} we combine the above two viewpoints with the general theory from \cite{DeC3} to construct a concrete implementation of the group von Neumann algebra of a linking quantum groupoid between $\widehat{SU_q(2)}$ and some locally compact quantum group $\mathbf{H}$.\\

\noindent In \emph{the fourth section}, we show that $\mathbf{H}$ is isomorphic to the dual of Woronowicz' $\widetilde{E}_q(2)$ quantum group. We end by showing that the linking quantum groupoid is in fact cleft, i.e.~ that it can be implemented by a unitary 2-cocycle for $\mathscr{L}^{\infty}(SU_q(2))$.\\

\noindent \qquad \emph{Conventions and remarks on notation}\\

\noindent For the rest of the paper, we fix a real number $0<q<1$. We then denote \[\lambda = (q-q^{-1})^{-1} < 0.\]

\noindent By $\iota$ we always mean the identity map. \\

\noindent The dual of a vector space $V$ will be denoted as $V^{\circ}$, since we preserve the symbol $*$ for anti-linear involutions.\\

\noindent When talking about Hopf algebras, we always assume that the base field is $\mathbb{C}$ and that the antipode is bijective. We will use the Sweedler notation for comultiplications and comodule structures (see \cite{Swe1}). We use \cite{Kli1} as our main reference to quantum groups.\\

\noindent We will denote quantum groups (and more general structures) by the symbols $\mathbb{G},\mathbb{H},\ldots$, and their duals by the symbols $\mathbf{G},\mathbf{H},\ldots$ (we give some more information concerning this notation in the first section). From the second section on, \emph{$\mathbb{G}$ will always denote the compact quantum group $SU_q(2)$}.\\

\noindent The following tensor product notations will be used. For algebraic tensor products between vector spaces, we use $\odot$. For minimal tensor products between C$^*$-algebras, or for tensor products between Hilbert spaces, we use $\otimes$. For spatial tensor products between von Neumann algebras, we use $\bar{\otimes}$. For tensor products between \emph{elements} of any these structures, we always use $\otimes$. We will also use the leg numbering notation as is customary in quantum group theory, but provide them with extra square brackets so as not to conflict with ordinary index notations (so if $x\in B(\mathscr{H}^{\otimes 2})$, we denote $x_{[23]}=1\otimes x \in B(\mathscr{H}^{\otimes 3})$ etc.).\\

\noindent When $\mathscr{H}$ is a Hilbert space, we write $B_0(\mathscr{H})$ for the space of compact operators. When $\xi,\eta\in \mathscr{H}$, we write $\omega_{\xi,\eta} \in B(\mathscr{H})_*$ for the normal functional \[\omega_{\xi,\eta}(x) = \langle x\xi,\eta\rangle,\] the scalar product being anti-linear in the second variable. Matrix units in $B(l^2(\mathbb{N}))$ are denoted as $e_{ij}$.\\

\noindent We will use the following notations from $q$-analysis. For $n\in \mathbb{N}\cup\{\infty\}$ and $a\in \mathbb{C}$, we denote \[(a;q)_{n} = \overset{n-1}{\underset{k=0}{\prod}} (1-q^ka),\] which determines analytic functions in the variable $a$ with no zeroes in the open unit disc. We write $E_{q}(z)$ for the formal power series \[E_q(z) = \sum_{k\in \mathbb{N}} \frac{q^{\frac{1}{2}k(k-1)}}{(q;q)_k} z^k.\] Then for any $a\in \mathbb{C}$, we have that $E_q(a)$ will be a convergent sum, with \[E_q(a) = (-a;q)_{\infty}.\] We write \[\left[ \begin{array}{c} n \\ m \end{array}\right]_{q} = \frac{(q;q)_n}{(q;q)_{n-m}(q;q)_m}\] for $m\leq n$ natural numbers.\\

\noindent We will need to use structures which are defined very similarly to each other. Then the names for the structures are often indexed, and when multiple structures are used together, we will index the objects associated to these structures with the corresponding index. However, when the structures appear isolated, we will refrain from indexing any of its associated structure. Also, when we index something with two indices which are the same, we will sometimes take the liberty of indexing with just one times this index symbol.\\

\section{Linking quantum groupoids and Morita equivalence of quantum groups}

\subsection{von Neumann algebraic quantum groups and linking quantum groupoids}

\begin{Def}(\cite{Kus2}) A von Neumann algebraic quantum group $(M,\Delta)$ consists of a von Neumann algebra $M$ and a normal unital $^*$-homomorphism $\Delta: M\rightarrow M\bar{\otimes} M$ satisfying the coassociativity condition \[(\Delta\otimes \iota)\Delta = (\iota\otimes \Delta)\Delta,\] and for which there exist nsf (normal semi-finite faithful) weights $\varphi$ and $\psi$ on $M$ such that for all normal states $\omega$ on $M$ and all $x\in M+$ we have \[\varphi((\omega\otimes \iota)\Delta(x)) = \varphi(x) \qquad \textrm{(left invariance)},\] \[\psi((\iota\otimes \omega)\Delta(x)) = \psi(x) \qquad \textrm{(right invariance)}.\]\end{Def}

\noindent We will follow the useful practice of writing a couple $(M,\Delta)$ as $\mathscr{L}^{\infty}(\mathbb{G})$, i.e.~ as if it concerns the space of bounded measurable functions on some `locally compact quantum group' $\mathbb{G}$. We then also freely use the notations $C_0^r(\mathbb{G})$ for the associated reduced C$^*$-algebraic quantum group (\cite{Kus1}), $\mathscr{L}(\mathbb{G})$ for the dual von Neumann algebraic quantum group $(\widehat{M},\widehat{\Delta}), \ldots$ In particular, we may write $\mathscr{L}^{\infty}(\mathbb{G})= \mathscr{L}(\widehat{\mathbb{G}})$ by the Pontryagin duality for locally compact quantum groups (\cite{Kus2}). \emph{For convenience of notation, we will denote the Pontryagin dual of $\mathbb{G}$ as $\mathbf{G}$ in stead of $\widehat{\mathbb{G}}$}, as to treat two quantum groups in duality on a more equal level qua notation. We will then further also talk about `a left corepresentation of $\mathscr{L}^{\infty}(\mathbb{G})$' as being `a left representation of $\mathbb{G}$', but will for example call `a left coaction of $\mathscr{L}^{\infty}(\mathbb{G})$ on a von Neumann algebra $N$' a `right action of $\mathbb{G}$ on $N$'. However, when also the algebra acted upon is interpreted as the space of functions on some `quantum space $\mathbb{X}$', we also call this `a left action of $\mathbb{G}$ on $\mathbb{X}$'. (We admit that this left-right terminology is a bit tedious, and we refrain from further comments on it in the future (concerning the terminology introduced in Theorem \ref{TheoProj} for example).)\\

\begin{Def}(\cite{Rie1}) Let $P$ and $M$ be two von Neumann algebras. A \emph{linking von Neumann algebra} between $P$ and $M$ consists of a von Neumann algebra $Q$ together with a self-adjoint projection $e\in Q$ and $^*$-isomorphisms $P\rightarrow eQe$ and $M\rightarrow (1-e)Q(1-e)$, such that both $e$ and $(1-e)$ are full projections (i.e. have central support equal to 1).\\

\noindent Two von Neumann algebras $P$ and $M$ are called \emph{W$^*$-Morita equivalent} if there exists a linking von Neumann algebra between them.
\end{Def}

\noindent We will mostly just identify $P$ and $M$ with their parts inside a linking von Neumann algebra, thus neglecting the identifying maps. We will also write $Q_{ij} = e_{ii}Qe_{jj}$ with $e_{11} = e$ and $e_{22}=1-e$, and \[ Q = \left(\begin{array}{ll} Q_{11} & Q_{12} \\ Q_{21} & Q_{22}\end{array}\right),\] in well-established and easily interpreted matrix notation. We will also talk simply of `a linking von Neumann algebra' (without specifying what the corners are) or of `a linking von Neumann algebra for the von Neumann algebra $M$' (without specifying the von Neumann algebra in the upper left corner). In fact, this terminology dictates the strongness of the isomorphism one is interested in (keeping none, one or both of the diagonals pointwise fixed).\\

\begin{Def}(\cite{DeC1}) Let $(P,\Delta_P)$ and $(M,\Delta_M)$ be von Neumann algebraic quantum groups. A \emph{von Neumann algebraic linking quantum groupoid} between $(P,\Delta_P)$ and $(M,\Delta_M)$ consists of a linking von Neumann algebra $(Q,e)$ between $P$ and $M$, together with a (non-unital) coassociative normal $^*$-homomorphism $\Delta_Q:Q\rightarrow Q\bar{\otimes} Q$ satisfying \[\Delta_Q(e) = e\otimes e,\qquad \Delta_Q(1-e)= (1-e)\otimes (1-e)\] and \[(Q_{11},(\Delta_Q)_{\mid Q_{11}})\cong (P,\Delta_P),\] \[(Q_{22},(\Delta_Q)_{\mid Q_{22}})\cong (M,\Delta_M)\] by the isomorphisms appearing in the definition of a linking von Neumann algebra. \end{Def}

\noindent We will then denote by $\Delta_{ij}: Q_{ij}\rightarrow Q_{ij}\bar{\otimes}Q_{ij}$ the restriction of $\Delta_Q$ to $Q_{ij}$. (One may view $Q_{ij}\bar{\otimes}Q_{ij}$ as $(e_{ii}\otimes e_{ii})(Q\bar{\otimes}Q)(e_{jj}\otimes e_{jj})$, this tensor product then being as well-behaved as the usual spatial tensor product between von Neumann algebras.) We also follow the same conventions as for linking von Neumann algebras, and will talk about `a von Neumann algebraic linking quantum groupoid' or `a von Neumann algebraic linking quantum groupoid for $(M,\Delta)$'.\\

\noindent It seems apt to interpret such a von Neumann algebraic linking quantum groupoid as the group von Neumann algebra of some `locally compact quantum groupoid $\mathbb{L}$' with two objects (see \cite{Bic1} for this interpretation in an algebraic (and dual) setting). We therefore also write von Neumann algebraic linking quantum groupoids in the form \[Q = \mathscr{L}(\mathbb{L}) = \left(\begin{array}{cc} \mathscr{L}(\mathbb{H}) & \mathscr{L}(\mathbb{X}) \\ \mathscr{L}(\overline{\mathbb{X}}) & \mathscr{L}(\mathbb{G})\end{array}\right),\] where for example $\mathscr{L}(\mathbb{X})$ should be seen as a weakly closed linear span of a space of arrows $\mathbb{X}$ from $\mathbb{G}$ to $\mathbb{H}$. The fact that $\mathbb{L}$ is a `linking groupoid' then says in a way that `any arrow (=element) in $\mathbb{G}$ can be written as the composition of an arrow in $\mathbb{X}$ and an arrow in $\overline{\mathbb{X}}$', and similarly for $\mathbb{H}$. Of course, for ordinary groups $G$ and $H$ this would be too strong a condition to be really interesting (since $G$ and $H$ would then automatically be isomorphic). It therefore seems better to interpret this situation classically in the more general settings of \emph{groupoids} $G$ and $H$, where the above condition precisely captures the notion of a `linking groupoid' and the associated `Morita equivalence' between groupoids. This is our motivation to call two locally compact quantum groups $\mathbb{G}$ and $\mathbb{H}$ \emph{Morita equivalent} when there exists a von Neumann algebraic linking quantum groupoid $\mathscr{L}(\mathbb{L})$ between $\mathscr{L}(\mathbb{G})$ and $\mathscr{L}(\mathbb{H})$. We then call $\mathbb{L}$ the associated \emph{locally compact (l.c.) linking quantum groupoid}. We could then also call the object $\mathbb{X}$ a \emph{locally compact quantum (bi-)torsor}. The usage of the predicate `locally compact' seems justified by the results in \cite{DeC1}, where it is shown that one also has associated C$^*$-algebraic descriptions (just as for the von Neumann algebraic quantum groups themselves). Finally, we note that the above objects can also be fitted into the framework of measured quantum groupoids (\cite{Les1},\cite{Eno1}), and could in fact be seen as the most easy kind of example of this structure, apart from the locally compact quantum groups themselves.\\

\noindent \emph{We will always use the above notation for linking quantum groupoids and its constituents. When changing the notation for $\mathbb{L}$ by using a different font or by adding an index, we make the same notational change for the constituents.}\\

\noindent The following notion was introduced in \cite{DeC3}.

\begin{Def} Let $\mathbb{G}$ be a locally compact quantum group. A \emph{(unitary) projective right representation} of $\mathbb{G}$ on a Hilbert space $\mathscr{H}$ is a right action \[\alpha: B(\mathscr{H})\rightarrow \mathscr{L}^{\infty}(\mathbb{G})\bar{\otimes} B(\mathscr{H})\] of $\mathbb{G}$ on $B(\mathscr{H})$.\end{Def}

\noindent This definition is motivated by the fact that if $G$ is a locally compact group and $\mathscr{H}$ a separable Hilbert space, there is a one-to-one correspondence between (anti-)homomorphisms $G\rightarrow \textrm{Aut}(B(\mathscr{H}))$ and (anti-)homomorphisms $G\rightarrow \mathcal{U}(\mathscr{H})/S^1$, the group of unitaries divided out by the circle group (imposing the appropriate continuity conditions). \\

\noindent The following theorem is one of the main results of \cite{DeC3} (Theorem 6.2, which uses a slightly different terminology), and its proof is quite involved. (An independent and easier proof, disregarding the weight structures, will appear in \cite{DeC5}).\\

\begin{Theorem}\label{TheoProj} Let $\mathbb{G}$ be a locally compact quantum group and $\alpha$ a projective right representation of $\mathbb{G}$ on a Hilbert space $\mathscr{H}$. Then there exists a l.c.~ linking quantum groupoid $\mathbf{L}$ for $\mathbf{G}$ and a unitary $\mathcal{G}\in \mathscr{L}(\mathbf{X})\bar{\otimes} B(\mathscr{H})$, such that \[(\Delta_{\mathscr{L}(\mathbf{X})}\otimes \iota)\mathcal{G} = \mathcal{G}_{[13]}\mathcal{G}_{[23]}\] and \[\alpha(x) = \mathcal{G}^*(1\otimes x)\mathcal{G} \qquad \textrm{for all }x\in B(\mathscr{H}).\] \end{Theorem}

\noindent We then call $\mathcal{G}$ an associated \emph{(unitary) left $\mathbf{X}$-representation for $\mathbb{G}$} (or a \emph{(unitary) left $\mathscr{L}(\mathbf{X})$-corepresentation for $\mathscr{L}^{\infty}(\mathbb{G})$}).\\

\noindent One should see the above Theorem as a generalization of the construction of (the cohomology class of) a 2-cocycle function $\Omega:G\times G\rightarrow S^1$ from a projective representation for a locally compact group $G$. In fact, the notion of a 2-cocycle still makes sense in the quantum setting.

\begin{Def}(\cite{Eno2}) A  unitary 2-cocycle for a von Neumann algebraic quantum group $(M,\Delta)$ is a unitary $\Omega\in M\bar{\otimes}M$ such that \[(\Omega\otimes 1)(\Delta\otimes \iota)(\Omega) = (1\otimes \Omega)(\iota\otimes \Delta)(\Omega).\]
\end{Def}

\noindent It turns out that, given such a 2-cocycle, one can then define a von Neumann algebraic linking quantum groupoid structure for $(M,\Delta)$ on $Q=M\otimes M_2(\mathbb{C})= \left(\begin{array}{ll} M & M \\M & M \end{array}\right)$ by putting \[\Delta_Q \left(\begin{array}{ll} w & x \\y & z \end{array}\right) =  \left(\begin{array}{rr} \Omega \Delta(w)\Omega^{*}  & \Omega \Delta(x) \\ \Delta(y) \Omega^{*} & \Delta (z) \end{array}\right),\] the existence of invariant weights on $(M,\Delta_{\Omega})$ (with $\Delta_{\Omega} = \Omega\Delta(\,\cdot\,)\Omega^{*}$) being the only non-trivial part to verify (see again \cite{DeC3}). In this way, they allow one to construct new locally compact quantum groups by a cocycle deformation or twisting (see e.g.~ \cite{Fim1} and \cite{Kas1} for the construction of some new locally compact quantum groups by this procedure; and see \cite{DeC2} and the present paper (Theorem \ref{TheoCoc}) for examples of such a twisting procedure turning the function algebra of a compact quantum group into the function algebra of a non-compact quantum group).\\

\noindent Conversely, linking quantum groupoids $\mathbb{L}$ whose associated linking von Neumann algebra is trivial (i.e.~ isomorphic to $\mathscr{L}(\mathbb{G})\otimes M_2(\mathbb{C})$) are easily seen to arise from a 2-cocycle for $\mathscr{L}(\mathbb{G})$, and we will call such von Neumann algebraic linking quantum groupoids \emph{cleft}. However, it is important to remark that the associated C$^*$algebraic Morita equivalences are in general \emph{not} trivial. We further remark that we do not know of any condition, \emph{directly in terms of a projective representation}, which guarantees that the associated von Neumann algebraic linking quantum groupoid will be cleft.\\

\noindent The following lemma will be useful in concretely identifying the group von Neumann algebra of a l.c.~ linking quantum groupoid associated to some given projective representation.

\begin{Lem}\label{LemProj} Let $\mathbb{G}$ be a locally compact quantum group, $\alpha$ a projective right representation of $\mathbb{G}$ on $l^{2}(\mathbb{N})$, and $\mathbf{L}$ and $\mathcal{G}$ resp.~ an associated l.c.~ linking quantum groupoid and $\mathbf{X}$-representation of $\mathbb{G}$. We write $x_n = (\iota\otimes \omega_{\delta_n,\delta_0})\mathcal{G} \in \mathscr{L}(\mathbf{X})$.\\

\noindent Let $(\widetilde{Q},f)$ be a linking von Neumann algebra for the von Neumann algebra $\mathscr{L}(\mathbf{G})$, and suppose that $(y_n)$ is a sequence of elements in $\widetilde{Q}_{12}$ such that $y_m^*y_n = x_m^*x_n$ for all $m,n\in \mathbb{N}$, and for which the set $\{y_n x\mid x \in M, n\in \mathbb{N}\}$ is $\sigma$-weakly dense in $\widetilde{Q}_{12}$. Then there exists a unique isomorphism \[\pi: (\widetilde{Q},f)\rightarrow (\mathscr{L}(\mathbf{L}),e)\] which identifies the canonical projections, pointwise fixes $\mathscr{L}(\mathbf{G})$, and such that $\pi(y_n)=x_n$ for all $n\in \mathbb{N}$. \end{Lem}

\begin{proof} It follows by standard von Neumann algebra techniques that there exists a unique (non-unital) normal $^*$-homomorphic embedding $\pi:\tilde{Q}\rightarrow \mathscr{L}(\mathbf{L})$ which pointwise fixes $\mathscr{L}(\mathbf{G})$ and such that $\pi(y_n)=x_n$. Then $\pi(f)\leq e$, and we should prove equality. For this it is sufficient to prove that the $\sigma$-weakly closed linear span $T$ of \[\{x_n x\mid x \in \mathscr{L}(\mathbf{G}), n\in \mathbb{N}\}\] equals $\mathscr{L}(\mathbf{X})$. But since $\mathcal{G}$ is a unitary, we will have $\sum_{n} x_nx_n^* =1_{\mathscr{L}(\mathbf{H})}=e$, where the infinite sum is taken in (say) the $\sigma$-weak topology. From this $T=\mathscr{L}(\mathbf{X})$ follows.\\
\end{proof}

\noindent We will also need the following density result, albeit in a simpler form.

\begin{Prop}\label{PropDens} Let $\mathbb{L}$ be a l.c.~ linking quantum groupoid for $\mathbb{G}$. Then the linear span of the set \[\{\Delta_{\mathscr{L}(\mathbb{X})}(x)(1\otimes y)\mid x\in \mathscr{L}(\mathbb{X}),y\in \mathscr{L}(\mathbb{G})\}\] is $\sigma$-weakly dense in $\mathscr{L}(\mathbb{X})\bar{\otimes} \mathscr{L}(\mathbb{X})$.\end{Prop}

\begin{proof} We first remark that the $\sigma$-weak closed linear span $T$ of the above set equals the one of the set \[\{\Delta_{\mathscr{L}(\mathbb{X})}(x)(y\otimes z)\mid x\in \mathscr{L}(\mathbb{X}),y\in \mathscr{L}(\mathbb{G})\},\] using that $\Delta_{\mathscr{L}(\mathbb{X})}(xy)=\Delta_{\mathscr{L}(\mathbb{X})}(x)\Delta_{\mathscr{L}(\mathbb{G})}(y)$ for $x\in \mathscr{L}(\mathbb{X})$ and $y\in \mathscr{L}(\mathbb{G})$, and the fact that the statement of the proposition holds when replacing $(\mathscr{L}(\mathbb{X}),\Delta_{\mathscr{L}(\mathbb{X})})$ by $(\mathscr{L}(\mathbb{G}),\Delta_{\mathscr{L}(\mathbb{G})})$ (see \cite{Kus2}, Proposition 1.4).\\

\noindent So we know now that $T$ is a $\sigma$-weakly closed right $\mathscr{L}(\mathbb{G})\bar{\otimes} \mathscr{L}(\mathbb{G})$-module. Since the $\sigma$-weakly closed span of $T^*T$ (resp.~ $TT^*$) contains the unit of $\mathscr{L}(\mathbb{G})\bar{\otimes} \mathscr{L}(\mathbb{G})$ (resp.~ $\mathscr{L}(\mathbb{H})\bar{\otimes} \mathscr{L}(\mathbb{H})$), the equality $T=\mathscr{L}(\mathbb{X})\bar{\otimes}\mathscr{L}(\mathbb{X})$ follows.

\end{proof}

\subsection{Co-linking weak Hopf algebras}

\begin{Def} A \emph{co-linking weak Hopf algebra} consists of four non-trivial unital algebras $H_{ij}$, $i,j\in\{0,1\}$, and eight unital homomorphisms $\Delta_{ij}^k: H_{ij}\rightarrow H_{ik}\odot H_{kj}$, $i,j,k\in\{0,1\}$, satisfying the coassociativity conditions \[(\Delta_{ik}^l\otimes \iota)\Delta_{ij}^k = (\iota\otimes \Delta_{lj}^k)\Delta_{ij}^l, \qquad i,j,k,l\in \{0,1\},\] and such that all maps \[H_{ij}\odot H_{kj} \rightarrow H_{ik}\odot H_{kj}: x\otimes y\rightarrow \Delta_{ij}^k(x)(1\otimes y),\]\[H_{ij}\odot H_{ik}\rightarrow H_{ik}\odot H_{kj}: x\otimes y\rightarrow \Delta_{ij}^k(x)(y\otimes 1)\] are bijective. (This will imply in particular that the $(H_{ii}^i,\Delta_{ii}^i)$ are Hopf algebras.)\\

\noindent If $(H_0,\Delta_0)$ and $(H_1,\Delta_1)$ are Hopf algebras, we say that the above co-linking weak Hopf algebra is \emph{a co-linking weak Hopf algebra between $H_0$ and $H_1$} when we are also given Hopf algebra isomorphisms $\Phi_0: (H_0,\Delta_0)  \rightarrow (H_{00},\Delta_{00}^0)$ and $\Phi_1:(H_1,\Delta_1)\rightarrow (H_{11},\Delta_{11}^1)$.\\

\noindent The substructure $(H_{ij},\Delta_{ij}^i,\Delta_{ij}^j)$ is called a \emph{bi-Galois object} between $(H_i,\Delta_i)$ and $(H_j,\Delta_j)$.\end{Def}

\noindent We call this a weak Hopf algebra because the separate pieces can be grouped together into \emph{a weak Hopf algebra with an invertible antipode} (\cite{Boh1}): one takes as the associated algebra $E$ the direct sum algebra $\oplus_{i,j\in\{0,1\}} H_{ij}$, and one defines the (non-unit preserving!) comultiplication $\Delta_E$ componentwise as \[\Delta_E(x_{ij})= \Delta_{ij}^0(x_{ij})+ \Delta_{ij}^1(x_{ij})\] for $x_{ij}\in H_{ij}$. One can in fact characterize independently the weak Hopf algebras which appear in this way (see \cite{DeC1}, Chapter 1 for details). They are then algebraic analogues of structures \emph{dual} to the von Neumann algebraic linking quantum groupoids of the previous subsection (i.e., one could picture them intuitively as a space of `polynomial functions' on a linking quantum groupoid). These dual structures were dealt with on an analytic level in \cite{DeC1} (see also \cite{DeC3}), but we will not be explicitly concerned with that side of the theory here. We also note that the above definition can be shown to be equivalent to that of a total Hopf-Galois system, as introduced in \cite{Gru1} (and extending work in \cite{Bic1}). We finally note that (bi-)Galois objects can also be defined independently (see \cite{Sch1}, and \cite{DeC1} for the correspondence with the above notion).\\

\begin{Not} When $(H_{ij},\Delta_{ij}^k)$ is a co-linking weak Hopf algebra, we will use the Sweedler-like notation \[ \Delta_{ij}^k(z) = z_{(1)ik}\otimes z_{(2)kj}, \qquad z\in H_{ij}\] for the coproduct.\end{Not}

\noindent We remark that the above notion has an obvious generalization to an `n$\times$n-co-linking weak Hopf algebra' consisting of $n^2$ algebras and $n^3$ comultiplications, or one large weak Hopf algebra with the $H_i$ at its `corners'. We also remark that one can consider co-linking weak Hopf $^*$-algebras and bi-$^*$-Galois objects, by putting a $^*$-structure on the algebras and imposing that it is respected by the comultiplication maps.\\

\begin{Prop}\label{PropAnt} When $(H_{ij},\Delta_{ij}^k)$ is a co-linking weak Hopf $^*$-algebra, there exist unique bijective anti-homomorphisms \[S_{ij}:H_{ij}\rightarrow H_{ji},\] each of which we call \emph{an antipode}, such that for all $i,j$ and $x\in H_{ii}$ we have \[S_{ij}(x_{(1)ij})x_{(2)ji} = \varepsilon_{ii}(x) = x_{(1)ij}S_{ji}(x_{(2)ji}).\]

\noindent They flip the comultiplication in the sense that \[\Delta_{ji}^k(S_{ij}(x)) = (S_{kj}\otimes S_{ik})(\Delta_{ij}^{k,\textrm{op}}(x)),\qquad \textrm{for all }x\in H_{ij}.\]

\noindent Furthermore, we have that for $x\in H_{ij}$, \[(S_{ji}(S_{ij}(x)^*))^* = x.\]\end{Prop}

\begin{proof} This can be pieced together from the results \cite{Sch2} and \cite{Boh1}. In any case, one can easily prove it by adapting the standard Hopf algebraic techniques.\end{proof}

\noindent There is also a counit $\varepsilon_E$ present on the weak Hopf algebra $E=\oplus H_{ij}$, given by \[x\rightarrow \delta_{ij}\varepsilon_i(x),\qquad x\in H_{ij}.\] One should however be careful with this notion, as the counit is not a homomorphism in the case of weak Hopf algebras.\\

\noindent Finally, there is the notion of an `adjoint action' for a co-linking weak Hopf algebra, whose technical term (referring to bi-Galois objects) is `the Miyashita-Ulbrich action'.

\begin{Def} Let $(H_{ij},\Delta_{ij}^k)$ be a co-linking weak Hopf algebra. The \emph{right Miyashita-Ulbrich action} of $(H_2,\Delta_2)$ on $H_{12}$ is the module algebra structure $\lhd$ on $H_{12}$ determined by \[x\lhd y = S_{21}(y_{(1)21})\cdot x\cdot y_{(2)12},\qquad x\in H_{12},y\in H_{22}.\]\end{Def}

\noindent In case there is a compatible $^*$-structure, we note that then \[x^* \lhd y = (x\lhd S(y)^*)^*,\] by the $^*$-identity for the antipode given in Proposition \ref{PropAnt}, i.e.~ $H_{12}$ is a module $^*$-algebra for $H_{22}$.

\subsection{On certain pairings between algebras}

\begin{Def} The \emph{$^*$-algebra of polynomial functions on $\mathbb{G}=SU_q(2)$} (\cite{Wor1}) is defined as the unital $^*$-algebra \[\textrm{Pol}(\mathbb{G})=\mathbb{C}\lbrack \mathbf{G}\rbrack=\textrm{Pol}_q(+,+)=\textrm{Pol}(SU_q(2)),\] generated (as a unital $^*$-algebra) by two generators $a_+$ and $b_+$ satisfying the relations \[\left\{\begin{array}{lllclll} a_+^*a_+ +\;\;\;\,b_+^*b_+ = 1 && \!\!\!\!\!\!\!\!\!\!\!\!a_+b_+ = q\;\;\;\,b_+a_+ \\a_+a_+^* + q^2 b_+b_+^* = 1 &&\!\!\!\!\!\!\!\!\!\!\!\! a_+^*b_+ = q^{-1}b_+a_+^* \\ &\!\!\!\!\!\!\!\!b_+b_+^* = b_+^*b_+.\end{array}\right.\]

\noindent We consider it as a Hopf $^*$-algebra by endowing it with the comultiplication map $\Delta_+$ satisfying \[\left\{\begin{array}{l} \Delta_+(a_+) = a_+\otimes a_+ + (- q b_+^*) \otimes b_+ \\ \Delta_+(b_+) = b_+\otimes a_+ + \;\;\;\;\quad a_+^*\otimes b_+.\end{array}\right.\]

\vspace{1cm}

\noindent The $^*$-algebra of \emph{polynomial functions on the quantum $\widetilde{E}(2)$ group} (\cite{Vak1},\cite{Cel1},\cite{Koe1}) is defined as the universal unital $^*$-algebra \[\textrm{Pol}(\mathbb{H})=\mathbb{C}\lbrack \mathbf{H}\rbrack = \textrm{Pol}_q(0,0)=\textrm{Pol}(\widetilde{E}_q(2))\] generated by elements $a_0$ and $b_0$, subject to the relations
\[\left\{\begin{array}{lllclll} a_0^*a_0 = 1 && \!\!\!\!\!\!\!\!\!\!\!\!a_0b_0 = q\;\;\;\,b_0a_0 \\a_0a_0^* = 1 &&\!\!\!\!\!\!\!\!\!\!\!\! a_0^*b_0 = q^{-1}b_0a_0^* \\ &\!\!\!\!\!\!\!\!b_0b_0^* = b_0^*b_0.\end{array}\right.\] We can make it into a Hopf $^*$-algebra  by endowing it with the comultiplication map $\Delta_0$ satisfying  \begin{eqnarray*} \Delta_0(a_0) &=& a_0\otimes a_0 \\ \Delta_0(b_0) &=& b_0\otimes a_0 + a_0^*\otimes b_0.\end{eqnarray*} \end{Def}

\noindent We note that \begin{eqnarray*} SU(2) &=& \{ \left(\begin{array}{rr} a & -\overline{b} \\ b & \overline{a} \end{array} \right) \mid |a|^2+|b|^2 = 1\},\\ \widetilde{E}(2)&\cong & \{ \left(\begin{array}{rr} a & 0 \\ b & \overline{a} \end{array} \right) \;\;\,\mid |a|^2\quad \quad\;\;= 1\},\end{eqnarray*}
\vspace{-0.2cm}
\noindent and
\vspace{-0.2cm}
\begin{eqnarray*} \!SU(1,1)\!&=&\!\{ \left(\begin{array}{rr} a & \overline{b} \\ b & \overline{a} \end{array} \right) \mid |a|^2-|b|^2 = 1\},\end{eqnarray*}
\noindent are precisely the Lie subgroups of $SL(2,\mathbb{C})$ keeping invariant the hermitian forms on $\mathbb{C}^2$ with symbol $(++)$, $(+0)$ and $(+-)$. This is one of the reasons for using the `sign notation' in the above and subsequent definitions.

\begin{Def} The \emph{quantized universal enveloping algebra of $su(2)$} (\cite{Kli1}) is defined as the universal algebra \[U_q(\mathbb{G}) = U_q(+,+) = U_q(su(2))\] generated by four elements $K_+,K_+^{-1},E_+$ and $F_+$ satisfying the commutation relations \[K_+^{-1}K_+ = 1 =  K_+K_+^{-1}, \;\;\;\; K_+E_+ = qE_+K_+,\;\;\;\;  K_+F_+ = q^{-1}F_+K_+,\] and \[\lbrack E_+,F_+\rbrack = \lambda(K_+^{2} - K_+^{-2}),\] where we recall that $\lambda = (q-q^{-1})^{-1}$.

\vspace{0.7cm}

\noindent The \emph{quantized universal enveloping algebra of the Lie algebra $e(2)$} (\cite{Vak1}) is defined as the universal algebra \[U_q(\mathbb{H}) = U_q(0,0) = U_q(e(2))\] generated by four elements $K_0,K_0^{-1},E_0$ and $F_0$ satisfying the commutation relations \[K_0^{-1}K_0 = 1 =  K_0K_0^{-1}, \;\;\;\; K_0E_0 = qE_0K_0,\;\;\;\;  K_0F_0 = q^{-1}F_0K_0,\] and \[\lbrack E_0,F_0\rbrack = 0.\]

\noindent Then both $U_q(+,+)$ and $U_q(0,0)$ can be made into Hopf $^*$-algebras by putting $K=K^*$, $E^*=F$ and \[\Delta(E)=E\otimes K+ K^{-1}\otimes E,\]\[\Delta(K)=K\otimes K,\] and then the antipode satisfies $S(K)=K^{-1}$, $S(E)=-qE$ and $S(F)=-q^{-1}F$.
\end{Def}

\begin{Prop}\label{PropPair} For $\mu\in \{0,+\}$, there is a pairing between the Hopf $^*$-algebras $U_q(\mu)$ and $\textrm{Pol}_q(\mu)$, uniquely determined by the formulas \[ \langle K,a \rangle = q^{-1/2}, \quad  \langle K,a^*\rangle = q^{1/2}, \quad \langle E,b \rangle = 1, \quad \langle F,(-qb^*)\rangle = 1,\] while all other possible pairings between generators are assigned zero.\end{Prop}

\begin{proof} For $\mu=+$, this is well-known (see \cite{Kli1}, Theorem 21). In any case, one can check it for both cases by using the argument at the beginning of section 4 of \cite{VDae2}.\end{proof}

\noindent To be clear, the fact that this is a pairing of Hopf $^*$-algebras (\cite{VDae2}) means that the associated maps \[\textrm{Pol}_q(\mu) \rightarrow U_q(\mu)^{\circ}\quad \textrm{and} \quad U_q(\mu)\rightarrow \textrm{Pol}_q(\mu)^\circ\] are unital homomorphisms (where the dual is given the convolution product), satisfying \[\langle x^*,y\rangle = \overline{\langle x,S(y)^*\rangle}, \quad \langle x,y^*\rangle = \overline{\langle S(x)^*,y\rangle}.\] Using an argument similar to the proof of Theorem 4.21 of \cite{Kli1}, one can show further that these pairings are non-degenerate. \\

\noindent Although we will not actually need the following proposition in its full generality, it provides a convenient reference for (easier) computations which otherwise would have to be treated separately.

\begin{Prop}\label{PropKli} For $n\in \mathbb{N}$, denote by $G_n(q)$ the number \[G_n(q) = \frac{(q^2;q^2)_n}{q^{\frac{1}{2}n(n-1)}(1-q^2)^n}.\]  Then we have for $s,r,t,n,l\geq 0$ and $m\in \mathbb{Z}$ that \[\langle K_{\mu}^mE_{\mu}^nF_{\mu}^l,a_{\mu}^rb_{\mu}^s(-qb_{\mu}^*)^t\rangle = \delta_{s,n}\delta_{t,l}q^{\frac{1}{2}m(-r+s-t)}q^{\frac{r}{2}(n+l)}G_n(q)G_l(q),\] in $\textrm{Pol}_q(\mu)$, $\mu\in\{0,+\}$. \\

\noindent Furthermore, for $\textrm{Pol}_q(0)$ this formula remains valid when $r\in \mathbb{Z}$. \end{Prop}

\begin{proof} We first argue that the above expression (for $r\in \mathbb{N}$) will be the same for both $\mu\in\{0,+\}$. Since the comultiplication of $U_q(\mu)$ restricts to the algebras generated by $K$,$K^{-1}$ and $E$, resp.~ $K,K^{-1}$ and $F$, and since the multiplication rules between $K$,$K^{-1}$ and $E$, resp.~ $K$,$K^{-1}$ and $F$ are the same in all cases, we only have to show that $\langle K^mE^nF^l,x\rangle$, with $x\in\{a,b,-qb^*\}$, is the same in all cases.\\

\noindent This is verified by direct computation, expanding the $x$ to $\Delta^{(m+n+l)}(x)$ and checking that only for one of the resulting terms, the pairing gives a non-zero number. Indeed, in case $x=a$, we see in this way that the above pairing is only non-zero in case $n=l=0$, in which case it equals $q^{-m/2}$. In case $x=b$, we get the number $q^{m/2}$, while for $x=-qb^*$, we get again $q^{-m/2}$.\\

\noindent The general formula in $(1)$ then follows from Proposition 22 in \cite{Kli1}, where the case $\mu=+$ is treated.\\

\noindent Since $E$ and $F$ commute in case $\mu=0$, the final statement is also immediate by applying the $^*$-operation (see also \cite{Koe1}, Proposition 3.2).\end{proof}

\section{The action on the standard Podle\'{s} sphere as a projective representation of $SU_q(2)$}

\subsection{Some facts on the standard Podle\'{s} sphere}

\begin{Def} (\cite{Wor1},\cite{Wor3}) We denote by $C(\mathbb{G})=C(SU_q(2))$ the (unique) C$^*$-algebraic quantum group associated to $\textrm{Pol}(SU_q(2))$, by $\mathscr{L}^{\infty}(\mathbb{G})=\mathscr{L}^{\infty}(SU_q(2))$ its associated von Neumann algebraic quantum group (containing then $\textrm{Pol}(SU_q(2))$ and $C(SU_q(2))$ as sub-$^*$-algebras), and by $\varphi$ the normal invariant state on the latter. \end{Def}

\begin{Def} The \emph{$^*$-algebra $\textrm{Pol}(S_{q0}^2)$ of polynomial functions on the standard Podle\'{s} sphere $S_{q0}^2$} is defined as the universal unital $^*$-algebra generated by elements $X$ and $Z$ satisfying the commutation relations \[ Z^*=Z, \;\;\;\;\;\; XZ = q^{2} ZX, \;\;\;\;\;\; X^*Z = q^{-2}ZX^*\] and \begin{eqnarray*} X^* X &=& Z- Z^2, \\ XX^* &=& q^2Z- q^4Z^2.\end{eqnarray*}

\noindent The \emph{C$^*$-algebra $C(S_{q0}^2)$ of continuous functions on the standard Podle\'{s} sphere $S_{q0}^2$} is defined as the universal enveloping C$^*$-algebra of $\textrm{Pol}(S_{q0}^2)$.\end{Def}

\noindent It is well-known that $C(S_{q0}^2)$ is isomorphic to $B_0(l^2(\mathbb{N})) + \mathbb{C}1$. A concrete identification is given by the unital $^*$-representation \[X\rightarrow \sum_{k\in \mathbb{N}_0} q^k \sqrt{1-q^{2k}}\; e_{k-1,k}\]\[Z\rightarrow \sum_{k\in \mathbb{N}} q^{2k}\, e_{kk},\] and we will then always regard $X$ and $Z$ as these concrete operators.\\

\noindent It is also well-known that $\textrm{Pol}(S_{q0}^2)\subseteq C(S_{q0}^2)$ can be imbedded into $\textrm{Pol}(SU_q(2))\subseteq C(SU_q(2))$ by the map \[X\rightarrow \tilde{X}=qb_+^*a_+,\]\[Z\rightarrow \tilde{Z}=b_+^*b_+.\] In fact, the range then becomes a (closed) left $^*$-coideal, which coincides with the $^*$-algebra of invariants for the action of the circle group $S^1$ on $\textrm{Pol}(SU_q(2))$ (resp.~ $C(SU_q(2))$) determined by $\theta_r(a_+)=e^{2\pi ir}a_+$, $\theta_r(b_+) = e^{2\pi ir}b_+$. We denote the corresponding continuous left action of $SU_q(2)$ on $S_{q0}^2$ by $\alpha$.\\

\noindent Finally, note that we get a state $\omega$ on $C(S_{q0}^2)$ by integrating out the action: \[\omega(x)1 = (\varphi\otimes \iota)\alpha(x).\] Then we can take the GNS-construction $\pi_{\omega}$ with respect to this state, and obtain the von Neumann algebra $\mathscr{L}^{\infty}(S_{q0}^2) = \pi_{\omega}(S_{q0}^2)''$, whose associated normal state we again denote as $\omega$. It is then easy to verify that $\mathscr{L}^{\infty}(S_{q0}^2)$ can be identified with $B(l^2(\mathbb{N}))$ in such a way as to extend the previous identification of $C(S_{q0}^2)$ with $B_0(\mathbb{N})+\mathbb{C}1$. The functional $\omega \in \mathscr{L}^{\infty}(S_{q0}^2)_*$ then coincides with $\textrm{Tr}(\,\cdot\, D)$, where $D$ is the trace class operator $(1-q^2)\textrm{Diag}(q^{2k})$ (one can prove this for example by using for example the techniques from the appendix of \cite{Wor3}). \\

\noindent The action $\alpha$ can then also be extended to an (ergodic) action \[\alpha:\mathscr{L}^{\infty}(S_{q0}^2)\rightarrow \mathscr{L}^{\infty}(SU_q(2))\bar{\otimes} \mathscr{L}^{\infty}(S_{q0}^2),\] which hence determines a unitary projective right representation of $SU_q(2)$. In the next sections, we find an explicit form for the associated von Neumann algebraic linking quantum groupoid \[\mathscr{L}(\mathbf{L})= \left(\begin{array}{cc} \mathscr{L}(\mathbf{H}) & \mathscr{L}(\mathbf{X}) \\\mathscr{L}(\overline{\mathbf{X}}) & \mathscr{L}(\mathbf{G})\end{array}\right),\] where we recall that $\mathbf{G}=\widehat{SU_q(2)}$. Note that by Proposition 10.3.4 of \cite{DeC1}, we already know that $\mathbf{H}$ will be a \emph{non-discrete} locally compact quantum group (and thus $\mathbb{H}$ not compact). We show in the final section that $\mathscr{L}(\mathbf{H})=\mathscr{L}^{\infty}(\mathbb{H})$ is indeed $\mathscr{L}^{\infty}(\widetilde{E}_q(2))$, justifying the notation used in the previous subsection.\\

\subsection{The standard Podle\'{s} sphere as a subquotient of a bi-Galois object}

\noindent The results in this subsection are dealt with in more detail in \cite{DeC4}.\\

\begin{Def} Let $\mu,\nu$ be two real numbers. We define $U_q(\mu,\nu)$ as the universal $^*$-algebra generated by four elements $K,K^{-1},E$ and $F$ satisfying the commutation relations \[K^{-1}K = 1 =  KK^{-1}, \;\;\;\; KE = qEK,\;\;\;\;  KF = q^{-1}FK,\]\[\;\;\;\; K^*=K,\;\;\;\; E^*=F\] and \[\lbrack E,F\rbrack = \lambda(\mu K^{2} - \nu K^{-2}).\] \end{Def}

\begin{Lem} There exists a family of unital $^*$-algebra homomorphisms \[\Delta_{\mu\nu}^\upsilon: U_q(\mu,\nu)\rightarrow U_q(\mu,\upsilon)\odot U_q(\upsilon,\nu)\] such that \[\Delta_{\mu\nu}^{\upsilon}(E_{\mu\nu}) = E_{\mu\upsilon}\otimes K_{\upsilon\nu}+K_{\mu\upsilon}^{-1}\otimes E_{\upsilon\nu},\] \[\Delta_{\mu\nu}^{\upsilon}(K_{\mu\nu}) = K_{\mu\upsilon}\otimes K_{\upsilon\nu}.\]
\end{Lem}
\begin{proof} By direct verification.\end{proof}

\noindent It is not difficult to see that for any $\mu>0$, we have $(U_q(\mu,\mu),\Delta_{\mu\mu}^{\mu})\cong (U_q(su(2)),\Delta)$, so we write $1=+$ again. We note further that for $\mu<0$ we will have $(U_q(\mu,\mu),\Delta_{\mu\mu}^{\mu})\cong (U_q(-),\Delta) = (U_q(su(1,1)),\Delta)$, the quantized universal enveloping algebra associated to $su(1,1)$ (\cite{Kli1}). But as we remarked in the introduction, we will treat this case in another paper. Note also that in the contraction procedure (\cite{Cel1}), one considers $U_q(0,0)$ as $\lim_{\mu\rightarrow 0} U_q(\mu,\mu)$.\\

\noindent To be more explicit, we have the respective Lie brackets \[\left.\begin{array}{lll} \lbrack E_{++},F_{++}\rbrack &=& \lambda(K_{++}^2-K_{++}^{-2}),\\ \lbrack E_{0+},F_{0+}\rbrack &=& \lambda(\qquad \;\,-K_{0+}^{-2}\,),\\ \lbrack E_{+0},F_{+0}\rbrack &=& \lambda (K_{+0}^2\;\;\qquad\;\;\;\;),\\ \lbrack E_{00},F_{00}\rbrack &=& 0,\end{array}\right.\] between the generators $E$ and $F$, which is another reason for using the sign notation.

\begin{Lem} The four-tuple $(U_q(\mu,\nu),\Delta_{\mu\nu}^\upsilon)$ with $\mu,\nu\in \{0,+\}$ is a co-linking weak Hopf $^*$-algebra.\end{Lem}

\noindent \emph{Remark:} The fact that $(U_q(0,+),\Delta_{0+}^0,\Delta_{0+}^+)$ is a bi-Galois object was first observed in \cite{Gun1}, Lemma 16. Since it preserves the $^*$-structure, it is automatically a bi-$^*$-Galois object.\\

\begin{proof} Define antipode maps $S_{\mu\nu}:U_q(\mu,\nu)\rightarrow U_q(\nu,\mu)$ by \begin{eqnarray*} S_{\mu\nu}(E_{\mu\nu}) &=& -qE_{\nu\mu},\\ S_{\mu\nu}(F_{\mu\nu})&=&-q^{-1}F_{\nu\mu}, \\ S_{\mu\nu}(K_{\mu\nu}) &=& K_{\nu\mu}^{-1}.\end{eqnarray*} Using the natural generalizations of the Hopf algebraic formulas concerning antipode and co-unit, one then easily constructs explicit inverses for the canonical maps appearing in the definition of a co-linking weak Hopf $^*$-algebra. See also \cite{DeC4}, Proposition 1.8.\end{proof}

\begin{Def} We define the \emph{Casimir element} of $U_q(\mu,\nu)$ to be the element \begin{eqnarray*} C_{\mu\nu} &=& EF + \lambda^2(q^{-1}\mu K^2+q \nu K^{-2})\\ &=& FE+\lambda^2(q\mu K^2+q^{-1}\nu K^{-2}).\end{eqnarray*}\end{Def}

\noindent As in the quantized universal enveloping algebra case, we have the following easy lemma.

\begin{Lem}\label{LemCas} The Casimir element $C_{\mu\nu}$ is a self-adjoint element in the center of $U_q(\mu,\nu)$. \end{Lem}

\begin{proof} The fact that $C_{\mu\nu}$ is self-adjoint is apparent on sight, while the statement that it lies in the center is immediate by checking that it commutes with $K$, $E$ and $F$.
\end{proof}

\noindent When $\tau$ is a real number, we will denote by $A_{\mu\nu}^\tau$ the quotient of $U_q(\mu,\nu)$ by the relation $C_{\mu\nu}= \tau$.\\

\noindent We note now that each $A_{\mu\nu}^\tau$ is endowed with more structure: since the right Miyashita-Ulbrich action by $U_q(\nu)$ on $U_q(\mu,\nu)$ is implemented by left and right multiplication maps, the ideal generated by the Casimir element $C_{\mu\nu}$ is invariant under it. Hence this action descends to a right $U_q(\nu)$-module $^*$-algebra structure on $A_{\mu\nu}^\tau$. \\

\noindent Take $\tau = q^{-1}\lambda^{2}$ (in fact, any non-zero number will work, but this value gives the nicest normalization), and then simply denote $A_{0+}^\tau$ as $A_{0+}$. The following is proven in \cite{DeC4} in the left setting. (We could also invite the reader to carry out the computations himself to experience the enjoyment of seeing terms cancel out.)

\begin{Prop}\label{PropIden2} There exists a faithful unital $^*$-homomorphic embedding of $\textrm{Pol}(S_{q0}^2)$ into $A_{0+}$ by $^*$-homomorphically extending the application \begin{eqnarray*} X &\rightarrow& \breve{X} = q^{1/2}(q^{-1}-q)E_{0+}K_{0+}^{-1},\\ Z&\rightarrow& \breve{Z} = K_{0+}^{-2},\end{eqnarray*} composed with the projection from $U_q(0,+)$ to $A_{0+}$. Moreover, endowing $\textrm{Pol}(S_{q0}^2)$ with the right $U_q(+)$-module structure \[y\lhd x = (\langle \,\cdot\,,x\rangle \otimes \iota)(\alpha(y)),\] with $\langle\,\cdot\,,\,\cdot\,\rangle$ the pairing from Proposition \ref{PropPair}, this is a right $U_q(+)$-module map.\end{Prop}

\section{Determining the quantum torsor $\mathbf{X}$.}

\subsection{The algebraic setting}

\noindent Denote $E=\underset{\mu,\nu\in\{0,+\}}{\oplus} U_q(\mu,\nu)$, with the $U_q(\mu,\nu)$ as defined in the previous section. We can then write the dual of $E$ in the form \[E^{\circ} = \left(\begin{array}{cc} U_q(0,0)^{\circ} & U_q(0,+)^{\circ} \\ U_q(+,0)^{\circ}& U_q(+,+)^{\circ}\end{array}\right),\] since we can endow it with the convolution multiplication dual to the comultiplication of $E$, and this will be of the form of a 2-by-2-matrix multiplication. More precisely, looking at the structure piecewise, we have the multiplication maps \[U_q(\mu,\upsilon)^{\circ}\odot U_q(\upsilon,\nu)^{\circ}\rightarrow U_q(\mu,\nu)^{\circ}: \omega_1\otimes \omega_2\rightarrow \omega_1\cdot\omega_2:\]\[ (\omega_1\cdot \omega_2)(x) = (\omega_1\otimes \omega_2)(\Delta_{\mu\nu}^{\upsilon}(x)),\qquad x\in U_q(\mu,\nu).\] We can then make it into a $^*$-algebra by the piecewise $^*$-operation \[*: U_q(\mu,\nu)^{\circ}\rightarrow U_q(\nu,\mu)^{\circ}:\] \[\omega^*(x) = \overline{\omega(S_{\nu\mu}(x)^*)},\qquad \omega\in U_q(\mu,\nu)^{\circ}, x\in U_q(\nu,\mu).\] \noindent Note in particular then that $U_q(\mu,\nu)$ is a $U_q(\mu)$-$U_q(\nu)$-bimodule.
\vspace{0.2cm}
\begin{Def} We define $\vartheta_{+0}$ (resp.~ $\vartheta_{0+}$) as the element in $U_q(+,0)^{\circ}$ (resp.~ $U_q(0,+)^{\circ}$) such that \[\langle \vartheta_{+0} , K_{+0}^mE_{+0}^nF_{+0}^l\rangle = \delta_{n,0}\delta_{l,0}\] (resp.~ $\langle \vartheta_{0+} , K_{0+}^mE_{0+}^nF_{0+}^l\rangle = \delta_{n,0}\delta_{l,0}$).
\end{Def}
\vspace{0.2cm}

\begin{Prop}\label{LemCom2} The following commutation relations hold:\[ a_{+}\vartheta_{+0} = \vartheta_{+0}a_{0},\quad b_+\vartheta_{+0} = \vartheta_{+0} b_0,\quad  b_+^*\vartheta_{+0} = \vartheta_{+0} b_0^*,\]\[a_{+}^*\vartheta_{+0} = \vartheta_{+0} (1-b_0^*b_0)a_0^*\\ =(1-b_+^*b_+)\vartheta_{+0}a_0^*.\]\end{Prop}

\begin{proof} It is immediate to see that for any $x\in \textrm{Pol}_q(+)$, we have \[\langle x\vartheta_{+0},K_{+0}^mE_{+0}^nF_{+0}^l\rangle = \langle x,K_{+}^mE_{+}^nF_{+}^l\rangle,\] and similarly for $x\in \textrm{Pol}_q(0)$ we have \[\langle \vartheta_{+0}x,K_{+0}^mE_{+0}^nF_{+0}^l\rangle = \langle x,K_{0}^mE_{0}^nF_{0}^l\rangle.\] Then the commutation relations follow immediately from Proposition \ref{PropKli} (where only the pairing with $a_+^*$ is not stated, but this is easily computed by hand or by referring to \cite{Kli1}, Proposition 22).\\
\end{proof}

\noindent We now remark that the multiplication of the co-linking weak Hopf algebra $E$ dualizes to `comultiplication maps' \[\Delta_{\mu\nu}: U_q(\mu,\nu)^{\circ}\rightarrow (U_q(\mu,\nu)\odot U_q(\mu,\nu))^{\circ}\] on the summands $U_q(\mu,\nu)^{\circ}$ of $E^{\circ}$. However, for the off-diagonal parts there does not seem to be a good `reduced sub-structure' to which it restricts to a genuine coalgebra structure (i.e.~ with the comultiplication ending up in an algebraic tensor product). In fact, we have the following proposition.

\begin{Prop}\label{PropIdenThet} The following identities hold, where $E_{q^{2}}$ is regarded as a power series, and with infinite sums converging in the weak topology: \begin{enumerate}
\item $\vartheta_{0+}\vartheta_{+0} = 1_{+}$, $\vartheta_{+0}\vartheta_{0+}=1_0$.
\item $\vartheta_{+0}^{*} = \vartheta_{0+}E_{q^{2}}(-q^{2}b_{+}^* b_{+})$,
\item $\Delta_{0+}(\vartheta_{+0}^*) = \sum_{p=0}^{\infty} (q^2;q^2)_p^{-1} ((a_0)^p\vartheta_{+0}^*b_+^p\otimes (a_0)^p\vartheta_{+0}^*(-qb_+^*)^p)$.
\item $(S_{+0}(\vartheta_{+0}))^{*}=S_{0+}(\vartheta_{+0}^*)=\vartheta_{+0}$.
\end{enumerate}
\end{Prop}

\noindent \emph{Remark:} The antipode $S_{+0}$ on $U_q(+,0)^{\circ}$ in the final expression is of course defined by the formula \[\langle S_{+0}(\omega),y\rangle = \langle \omega,S_{0+}(y)\rangle, \qquad \omega\in U_q(+,0)^{\circ},y\in U_q(0,+),\]
and similarly for $S_{0+}$.

\begin{proof} \noindent The first identity is immediate.\\

\noindent We now compute first the final identities. We have \begin{eqnarray*}\langle  (S_{+0}(\vartheta_{+0}))^*,K_{+0}^mE_{+0}^nF_{+0}^l\rangle &=&  \langle \vartheta_{+0}, (K_{+0}^mE_{+0}^nF_{+0}^l)^*\rangle \\ &=& \langle \vartheta_{+0},E_{+0}^lF_{+0}^nK_{+0}^m\rangle \\ &=& \delta_{n,0}\delta_{l,0}\\ &=& \langle \vartheta_{+0},K_{+0}^mE_{+0}^nF_{+0}^l\rangle,\end{eqnarray*} and similarly for the other one.\\

\noindent For the remaining two identities, we first remark that \[\langle \vartheta_{+0},F_{+0}^lE_{+0}^n\rangle = \delta_{nl}\, q^n\,\frac{(q^2;q^2)_n}{(1-q^2)^{2n}},\] which is easily computed by induction (see also Lemma 3.2.2.(ii) in \cite{Maj1}). Using the full formula for the pairing from Proposition \ref{PropKli}, the second and third identity are then straightforwardly computed by taking the pairing with an element $K_{0+}^mE_{0+}^nF_{0+}^l$ or $K_{0+}^mE_{0+}^nF_{0+}^l\otimes K_{0+}^{m'}E_{0+}^{n'}F_{0+}^{l'}$ (first transporting $b_+$ and $b_+^*$ across $\theta_{+0}^*$ into $b_0$ and $b_0^*$ in the second expression, which is allowed by Proposition \ref{LemCom2}).
\end{proof}

\noindent \emph{Remark:} The element $\Delta_{+0}(\vartheta_{+0})$ is an (invertible) algebraic 2-cocycle functional $\omega$ associated the bi-Galois object $U_q(+,0)$ (see \cite{Kas1} and \cite{Aub1}). But we see from the above Proposition that it is not \emph{unitary}. However, it does satisfy the formula $\omega(x^*,y^*) = \overline{\omega(y,x)}$, i.e.~ is \emph{real} in the terminology of \cite{Maj1}, Proposition 2.3.7. We will see in the next subsection that in the analytic picture, $\omega$ can be implemented by an operator, but \emph{this operator will no longer be invertible}. Since the 2-cocycle twisting of von Neumann algebraic quantum groups has only been studied in full generality for unitary 2-cocycles, we thus have resorted to the more general theory of Morita equivalence of locally compact quantum groups, even though, \emph{in the end}, our example \emph{will} be of the form of a unitary 2-cocycle twist (but in a very non-natural way). This final result then leads us again to a question concerning the 2-cocycle functional $\omega$: will it be coboundary equivalent to a \emph{unitary} 2-cocycle functional? Given the non-naturality of our operator algebraic 2-cocycle, this is not at all clear.\\

\noindent In the next subsection, we will also need the following easy result. We first introduce a definition/notation.

\begin{Def} We define $\mathbb{C}\lbrack \mathbf{X}\rbrack $ as the linear span of the set \[\{a_0^r\vartheta_{+0}^*b_+^s(b_+^*)^t\mid r\in \mathbb{Z},s,t\in \mathbb{N}\}.\]\end{Def}

\begin{Lem}\label{LemBas} The vector space $\mathbb{C}\lbrack \mathbf{X}\rbrack$ is a right $\mathbb{C}\lbrack \mathbf{G}\rbrack$-module with the set \[\{a_0^r\vartheta_{+0}^*b_+^s(b_+^*)^t\mid r\in \mathbb{Z},s,t\in \mathbb{N}\}\] as a vector space basis.\end{Lem}

\begin{proof} Using the commutation relations in Proposition \ref{LemCom2}, one checks that $\mathbb{C}\lbrack \mathbf{X}\rbrack$ is globally invariant under right multiplication with the generators of $\mathbb{C}\lbrack \mathbf{G}\rbrack = \textrm{Pol}(+)$. Hence it is a right $\mathbb{C}\lbrack \mathbf{G}\rbrack$-module.\\

\noindent For the linear independence statement, we will give a proof which can be recycled later on. Suppose that $\sum_{r,s,t} c_{r,s,t} a_0^r\vartheta_{+0}^*b_+^s(b_+^*)^t =0$, where all but a finite number of $c_{r,s,t}$ are zero. Multiplying to the left with a high enough power of $a_0^*$, we may assume that $c_{r,s,t}=0$ for $r> 0$. Then by Lemma \ref{LemCom2}, we have $\sum_{r,s,t\in \mathbb{N}} c_{-r,s,t} \vartheta_{+0}^*(a_+^*)^rb_+^s(b_+^*)^t =0$. Multiplying to the left with $\vartheta_{0+}^*$, we have $\sum_{r,s,t} c_{-r,s,t} (a_+^*)^rb_+^s(b_+^*)^t =0$ by the first identity in Proposition \ref{PropIdenThet}.  Since it is well-known that the $(a_+^*)^rb_+^s(b_+^*)^t$ are linearly independent, this proves that all $c_{r,s,t}$ are zero.

\end{proof}

\noindent We will need the following computation in the next part.

\begin{Lem}\label{LemComp1} Let $m\in \mathbb{Z}$ and $s,n,l\in\mathbb{N}$. Then  \[\langle a_0^s \vartheta_{+0}^*b_+^s, K_{0+}^mF_{0+}^lE_{0+}^n\rangle =   \delta_{l,0} \delta_{s,n} q^{n/2}  (q^2;q^2)_n (1-q^2)^{-n}.\] \end{Lem}

\begin{proof} First moving the $^*$ to the other side as $S(\,\cdot\,)^*$, this again follows immediately from Lemma \ref{LemCom2} (moving $b_+$ across $\vartheta_{+0}$ into $b_0$, and then deleting $\vartheta_{+0}$ as to obtain a pairing between elements of $\textrm{Pol}(0)$ and $U_q(0)$), and then applying the formula in Proposition \ref{PropKli}.
\end{proof}

\subsection{The analytic setting}

\noindent It is well-known that $\mathscr{L}^{\infty}(SU_q(2))$ can be represented faithfully on $l^2(\mathbb{N})\otimes l^2(\mathbb{Z})$ by the representation \[a_+\rightarrow (\sum_{k\in \mathbb{N}_0} \sqrt{1-q^{2k}}e_{k-1,k})\otimes 1,\]\[b_+\rightarrow (\sum_{k\in \mathbb{N}} q^ke_{kk})\otimes S,\] with $S$ denoting the bilateral shift forward. We will in the following always identify $a_+$ and $b_+$ with their images in this representation. We then further extend our sign notation.

\begin{Not} We denote \[\mathscr{H}_{++} = l^2(\mathbb{N})\otimes l^2(\mathbb{Z})\] and \[\mathscr{H}_{0+} = l^2(\mathbb{Z})\otimes l^2(\mathbb{Z}).\] \end{Not}

\begin{Def} We define $v_0$ as the unitary operator $S^*\otimes 1$ on $\mathscr{H}_{0+}$, with $S^*$ the backward bilateral shift.\\

\noindent We define $L_{0+}$ as the bounded operator $\mathscr{H}_+\rightarrow \mathscr{H}_{0+}$ determined by \[L_{0+}(e_n^{+}\otimes e_k^+) = \left(\frac{(q^2;q^2)_{\infty}}{(q^2;q^2)_n}\right)^{1/2} e_n^{0+}\otimes e_k^{0+}.\]\end{Def}

\noindent More succinctly, $L_{0+} = uE_{q^2}(-q^2b_+^*b_+)^{1/2} = u(q^2b_+^*b_+;q^2)_{\infty}^{1/2}$, with $u$ the isometric imbedding \[u:\mathscr{H}^+\rightarrow \mathscr{H}^{0+}:e_n^+\otimes e_k^+ \rightarrow e_n^{0+}\otimes e_k^{0+}.\] We remark that this operator (or rather $L_{0+}u^*$) also appears in \cite{Wor6}, as the operator $t^{1/2}$ (formula (28)).

\begin{Lem}\label{LemCom1} The following commutation relations hold: \[L_{0+}a_+^* = v_0^*L_{0+},\] \[L_{0+}a_+ = v_0L_{0+}(1-b_+^*b_+).\]\end{Lem}

\begin{proof} This follows by an easy computation, which we leave to the reader.\end{proof}

\noindent By induction, we then get that \[L_{0+}a_+^n = v_0^nL_{0+}(b_+^*b_+;q^{-2})_{n}.\] Note however that we cannot move $v_0$ across $L_{0+}$, as $v_0$ shifts the range of $L_{0+}$ downwards. Also note that $1-b_+^*b_+$ has the same kernel as $a_+$, which is part of what makes the above commutation relation work.

\begin{Lem} The linear span of $\{v_0^rL_{0+}b_+^s(b_+^*)^t\mid r\in \mathbb{Z},s,t\in \mathbb{N}\}$ is a right $\textrm{Pol}(SU_q(2))$-module, isomorphic to $\mathbb{C}\lbrack \mathbf{X}\rbrack$ by the map \[\pi: a_0^r\vartheta_{+0}^*b_+^s(b_+^*)^t\rightarrow v_0^rL_{0+}b_+^s(b_+^*)^t.\]\end{Lem}

\begin{proof} By the commutation relations of Lemma \ref{LemCom1}, and because $L_{0+}^*L_{0+}$ is invertible, we can simply copy the proof of Lemma \ref{LemBas}, replacing $\vartheta_{+0}^*$ by $L_{0+}$ and $a_0$ by $v_0$, to get that $\{v_0^rL_{0+}b_+^s(b_+^*)^t\mid r\in \mathbb{Z},s,t\in \mathbb{N}\}$ is a basis for its linear span, and to see that this linear span is a right $\textrm{Pol}(SU_q(2))$-module by right multiplication. By the same commutation relations, the existence of the isomorphism $\pi$ is then also immediate.
\end{proof}

\noindent In the following, we will drop the symbol $\pi$ for the isomorphism, and thus identify $\mathbb{C}\lbrack \mathbf{X}\rbrack$ with a space of operators $\mathscr{H}_{+}\rightarrow \mathscr{H}_{0+}$. However, we keep using the distinct symbols $L_{0+}$, $\vartheta_{0+}$, $\ldots$ to clarify whether we consider something as an operator or a functional.\\

\noindent Define now $\widetilde{Q}$ as the linking von Neumann algebra \[\left(\begin{array}{cc} B(l^2(\mathbb{Z})) &  B(l^2(\mathbb{N}),l^2(\mathbb{Z}))\\  B(l^2(\mathbb{Z}),l^2(\mathbb{N})) & B(l^2(\mathbb{N})) \end{array}\right)\bar{\otimes} \mathscr{L}(\mathbb{Z}),\] represented on $\left(\begin{array}{c} \mathscr{H}_{0+} \\ \mathscr{H}_+\end{array}\right)$. We may identify $\widetilde{Q}_{22}$ with $\mathscr{L}^{\infty}(SU_q(2))$. Note further then that $\mathbb{C}\lbrack \mathbf{X}\rbrack$ is $\sigma$-weakly dense in $\widetilde{Q}_{12}$.\\

\begin{Prop}\label{PropRow} Define \[y_n = (q^2;q^2)_n^{-1/2} v_0^n L_{0+}b_+^n \in \widetilde{Q}_{12}.\] Then this sequence satisfies the conditions of Lemma \ref{LemProj} with respect to the coaction $\alpha$ of $\mathscr{L}^{\infty}(SU_q(2))$ on $B(l^2(\mathbb{N}))\cong \mathscr{L}^{\infty}(S_{q0}^2)$. \end{Prop}

\begin{proof} The density condition in that lemma is easily verified to hold true. We only have to prove then the identities $x_m^*x_n = y_m^*y_n$.\\

\noindent First observe that, since $\alpha(x)=\mathcal{G}^*(1\otimes x)\mathcal{G}$ for $x\in L^{\infty}(S_{0q}^2)$, we will have \[\alpha(e_{00}) = \sum_{m,n} x_m^*x_n \otimes e_{mn} \] as a $\sigma$-weakly converging sum, so that \[x_m^*x_n = (\iota\otimes \omega_{\delta_n,\delta_m})\alpha(e_{00}).\]

\noindent Now $e_{00} = \lim_{n\rightarrow \infty} Z^n$ in norm. Taking the copy $\tilde{Z}=b_+^*b_+$ of $Z$ in $\mathscr{L}^{\infty}(SU_q(2))$, we can write \begin{eqnarray*} \alpha(Z)^n &\cong & \Delta(\tilde{Z})^n \\ &=& \Delta(b_+^*)^n\Delta(b_+)^n \\ &=& \sum_{m,k=0}^{n} \qbin{n}{m}  \qbin{n}{k} (b_+^*)^{n-k} a_+^k(a_+^*)^m b_+^{n-m} \otimes (a_+^*)^{n-k}(b_+^*)^k b_+^ma_+^{n-m}, \end{eqnarray*} where we have used the well-known formula for the $n$-th power of $q^2$-commuting variables. \\

\noindent Now we remark that in fact $C(SU_q(2))\subseteq \mathscr{L}^{\infty}(S_{0q}^2)\otimes C(S^1) \cong C(S^1,\mathscr{L}^{\infty}(S_{0q}^2))$, and that we have a surjective $^*$-homomorphism $\theta$ from $C(SU_q(2))$ to $\mathscr{L}^{\infty}(S_{0q}^2)$ by evaluating in zero. Denoting \begin{eqnarray*} A&=&\sum_{k\in \mathbb{N}_0} \sqrt{1-q^{2k}}e_{k-1,k},\\ B&=& Z^{1/2}\end{eqnarray*} as elements in $\mathscr{L}^{\infty}(S_{0q}^2)$, we have $\theta(b_+) = B$ and $\theta(a_+)=A$, so that $\theta(\tilde{Z}) = Z$ and $\theta(\tilde{X})=X$. So \[\alpha(Z)^n = \sum_{m,k=0}^{n} \qbin{n}{m}  \qbin{n}{k} (b_+^*)^{n-k} a_+^k(a_+^*)^m b_+^{n-m} \otimes (A^*)^{n-k}B^kB^mA^{n-m},\] and then, for $s,t\in \mathbb{N}$ and $n\geq s,t$, \begin{eqnarray*} (\iota\otimes \omega_{\delta_t,\delta_s})(\alpha(Z)^n) &=& \sum_{m=n-t}^{n} \sum_{k=n-s}^n \frac{(q^2;q^2)_t^{1/2}}{(q^2;q^2)_{t-(n-m)}^{1/2}} \cdot \frac{(q^2;q^2)_s^{1/2}}{(q^2;q ^2)_{s-(n-k)}^{1/2}} q^{m(t-(n-m))} q^{k(s-(n-k))} \cdot\\ && \qquad \qquad \delta_{t-(n-m),s-(n-k)}\qbin{n}{m}\qbin{n}{k} (b_+^*)^{n-k}a_+^k(a_+^*)^mb_+^{n-m} \\ &=& \sum_{r=0}^{\textrm{min}\{s,t\}} \frac{(q^2;q^2)_t^{1/2}}{(q^2;q^2)_{r}^{1/2}} \cdot \frac{(q^2;q^2)_s^{1/2}}{(q^2;q^2)_{r}^{1/2}} q^{r(n-(t-r))}q^{r(n-(s-r))}\cdot \\ && \qquad  \qquad \qquad \qbin{n}{n-(t-r)}\qbin{n}{n-(s-r)} \\ && \qquad \qquad \qquad \qquad \qquad (b_+^*)^{s-r}a_+^{n-(s-r)}(a_+^*)^{n-(t-r)}b_+^{t-r}.\end{eqnarray*}

\noindent Taking the limit $n\rightarrow \infty$, we see that only the term $r=0$ survives, and we finally get \[x_s^*x_t = (q^2;q^2)_t^{-1/2}(q^2;q^2)_s^{-1/2} (b_+^*)^s (\lim_{n\rightarrow \infty} a_+^{n-s}(a_+^*)^{n-t})b_+^t.\]

\noindent On the other hand, by Lemma \ref{LemCom1}, \begin{eqnarray*} y_s^*y_t &=& (q^2;q^2)_s^{-1/2} (q^2;q^2)_t^{-1/2} (b_+^*)^sL_{0+}^*(v_0^*)^s v_0^t L_{0+}b_+^t \\ &=&  (q^2;q^2)_s^{-1/2} (q^2;q^2)_t^{-1/2} (b_+^*)^sa_+^t L_{0+}^* L_{0+}(a_+^*)^s b_+^t \\ &=&  (q^2;q^2)_s^{-1/2} (q^2;q^2)_t^{-1/2} (b_+^*)^sa_+^t(q^2b_+^*b_+;q^2)_{\infty}(a_+^*)^s b_+^t.\end{eqnarray*} An easy computation reveals then that \[ \lim_{n\rightarrow \infty} a_+^{n-s}(a_+^*)^{n-t} = a_+^t(q^2b_+^*b_+;q^2)_{\infty}(a_+^*)^s,\] which finishes the proof.

\end{proof}

\noindent \emph{Remark:} The value of $\alpha(e_{00})$ could also have been derived directly from the formula for $\Delta(e_{00})$ in \cite{Lan1} (using Theorem 4.1 and formula (4.7)). However, our proof is somewhat easier, and avoids for example transformation laws for basic hypergeometric series, instead using only very elementary $q$-analysis.\\

\noindent By the previous proposition and Lemma \ref{LemProj}, we may thus identify $\widetilde{Q}$ with $\mathscr{L}(\mathbf{L})$ in the foregoing manner, and we will drop the distinction from now on.\\

\noindent We can now make a $^*$-representation $\Theta$ of $U_q(0,+)$ on the pre-Hilbert space $\mathbb{C}\lbrack\mathbb{N}\rbrack\subseteq l^2(\mathbb{N})$ consisting of finite linear combinations of the basis vectors, uniquely determined by the properties that $\Theta(K_{0+})$ is a positive operator and such that, using the notation of Proposition \ref{PropIden2}, \[\Theta(\breve{X})=X_{\mid \mathbb{C}\lbrack \mathbb{N}\rbrack}\quad \textrm{and}\quad \Theta(\breve{Z})=Z_{\mid \mathbb{C}\lbrack \mathbb{N}\rbrack}.\] To wit, the image of $K_{0+}$ thus becomes the diagonal operator for which $\Theta(K_{0+})e_n = q^{-n}e_n$.\\

\noindent We recall that we denote by $\mathcal{G}$ the $\mathbf{X}$-representation of $SU_q(2)$ associated to the projective representation $\alpha$ (see Theorem \ref{TheoProj}). In the following proposition, we write $\mathcal{G}_{r,s}\in \mathscr{L}(\mathbf{X})$ for the element $(\iota\otimes \omega_{\delta_s,\delta_r})(\mathcal{G})$, so then \[\mathcal{G} = \sum_{r,s\in\mathbb{N}} \mathcal{G}_{r,s}\otimes e_{rs}\] in the $\sigma$-weak topology on $\mathscr{L}(\mathbf{X})\bar{\otimes} B(l^2(\mathbb{N}))$.

\begin{Prop}\label{PropVal} For all $r,s\in\mathbb{N}$, we have that $\mathcal{G}_{r,s}\in \mathbb{C}\lbrack \mathbf{X}\rbrack$, and \[\langle \mathcal{G}_{r,s},y\rangle = \Theta(y)_{r,s}\qquad \textrm{for all }y\in U_q(0,+).\]\end{Prop}

\begin{proof}

\noindent We give a proof by induction on $r$.\\

\noindent By Proposition \ref{PropRow}, the element $\mathcal{G}_{0,s}$ equals $(q^2;q^2)_s^{-1/2} v_0^sL_{0+}b_+^s$, which belongs to $\mathbb{C}\lbrack \mathbf{X}\rbrack$. So we should see if \[\langle a_0^s \vartheta_{+0}^*b_+^s, y\rangle = (q^2;q^2)_s^{1/2} \Theta(y)_{0,s}\] for $y\in U_q(0,+)$.\\

\noindent Choose $y$ of the form $K_{0+}^mF_{0+}^lE_{0+}^n$ for $m\in \mathbb{Z},n,l\in\mathbb{N}$. Then \[\langle a_0^s \vartheta_{+0}^*b_+^s, K_{0+}^mF_{0+}^lE_{0+}^n\rangle = \delta_{l,0} \delta_{s,n} q^{n/2}  (q^2;q^2)_n (1-q^2)^{-n}\] by Lemma \ref{LemComp1}.\\

\noindent On the other hand, we have that \begin{eqnarray*} \Theta(K_{0+}^mF_{0+}^lE_{0+}^n)_{0,s} &=& \omega_{\delta_s,\delta_0} (Z^{-\frac{m}{2}} (-q^{-1/2}\lambda Z^{-1/2}X^*)^l (-q^{-1/2}\lambda X Z^{-1/2})^n) \\ &=& \delta_{l,0} \delta_{s,n}  (-1)^nq^{-n/2} \lambda^n q^{-\frac{n(n+1)}{2}} \omega_{\delta_n,\delta_0}(X^n) \\ &=& \delta_{l,0} \delta_{s,n} q^{-n/2} (q^{-1}-q)^{-n} q^{-\frac{n(n+1)}{2}}    \omega_{\delta_n,\delta_0}(X^n) \\ &=& \delta_{l,0} \delta_{s,n} q^{-n/2} (q^{-1}-q)^{-n} q^{-\frac{n(n+1)}{2}}  q^{\frac{n(n+1)}{2}} (q^2;q^2)_n^{1/2}\\ &=&  \delta_{l,0} \delta_{s,n} q^{n/2} (1-q^2)^{-n} (q^2;q^2)_n^{1/2}.\end{eqnarray*}

\noindent This proves the case $r=0$.\\

\noindent Now suppose that the condition and formula hold for some $r$ and all $s$. Then from the fact that \[(1\otimes X)\mathcal{G} = \mathcal{G}\alpha(X),\] we conclude that \[ \mathcal{G}_{r+1,t} = q^{-(r+1)}(1-q^{2(r+1)})^{-1/2}\sum_{s} \omega_{\delta_t,\delta_s}(X_{(0)})\;\mathcal{G}_{r,s}X_{(-1)},\] where we remark that there are only finitely many non-zero terms in the right hand side. From this, we already see that $\mathcal{G}_{r+1,t}$ lies in $\mathbb{C}\lbrack \mathbf{X}\rbrack$. Using Proposition \ref{PropIden2}, we further compute that for $x\in U_q(0,+)$, we have \begin{eqnarray*} \langle \sum_{s} \omega_{\delta_t,\delta_s}(X_{(0)})\;\mathcal{G}_{r,s}X_{(-1)},x\rangle &=& \sum_{s} \omega_{\delta_t,\delta_s}(X_{(0)}) \;\langle \mathcal{G}_{r,s},x_{(1)0+}\rangle \;\langle X_{(-1)},x_{(2)+}\rangle \\ &=&  \sum_{s}  \omega_{\delta_t,\delta_s}(X\lhd x_{(2)+})\; \Theta(x_{(1)0+})_{r,s}  \\ &=&   \sum_{s}  \Theta(\breve{X}\lhd x_{(2)+})_{s,t}\; \Theta(x_{(1)0+})_{r,s}  \\ &=& \Theta(x_{(1)0+}S_{+0}(x_{(2)+0})\,\breve{X}\,x_{(3)0+})_{r,t} \\ &=& \Theta(\breve{X}x)_{r,t} \\ &=& (X\Theta(x))_{r,t} \\ &=& q^{r+1} (1-q^{2(r+1)})^{1/2} \Theta(x)_{r+1,t}.\end{eqnarray*}

\noindent This concludes the proof.

\end{proof}

\noindent We can then use this proposition to explicitly compute the comultiplication $\Delta_{0+}$ on $\mathscr{L}(\mathbf{X})$.

\begin{Theorem}\label{TheoComu} The comultiplication $\Delta_{0+}$ on $\mathscr{L}(\mathbf{X})$ is determined on $\mathbb{C}\lbrack \mathbf{X}\rbrack$ by the formula \[\Delta_{0+}(v_0^rL_{0+}) = (v_0^r\otimes v_0^r)\cdot (\sum_{p=0}^{\infty} (q^2;q^2)_p^{-1} \;v_0^pL_{0+}b_+^p\otimes v_0^pL_{0+}(-qb_+^*)^p),\] where the sum in the right hand side converges in norm.
\end{Theorem}

\begin{proof} \noindent The norm-convergence of the right hand side expression is of course immediate.\\

\noindent We first prove the identity part of the theorem in the case $r=0$. We make the following computation: for $x\in U_q(0,+)$, and considering again the $^*$-representation $\Theta$ introduced before Proposition \ref{PropVal}, we have, by that same proposition, \begin{eqnarray*} \Theta(x)_{r,0} &=& \overline{\Theta(x^*)_{0,r}} \\ &=& (q^2;q^2)_r^{-1/2}\langle a_0^r\vartheta_{+0}^*b_+^r,x^*\rangle \\ &=& (q^2;q^2)_r^{-1/2}\langle (S_{0+}(a_0^r\vartheta_{0+}^*b_+^r))^*,x\rangle \\ &=& (q^2;q^2)_r^{-1/2} (-q)^r\langle a_0^r\vartheta_{0+}^* (b_+^*)^r ,x\rangle, \end{eqnarray*}where we used the final identity in \ref{PropIdenThet} for the last step. Hence \[\mathcal{G}_{r,0} = (q^2;q^2)_r^{-1/2} v_0^rL_{+0}(-qb_+^*)^r.\] Since $L_{0+}= \mathcal{G}_{00}$ and $\Delta_{0+}(\mathcal{G}_{00}) = \sum_{p\in \mathbb{N}} \mathcal{G}_{0p}\otimes \mathcal{G}_{p0}$ in the $\sigma$-weak topology, we have proven the stated identity in case $r=0$.\\

\noindent Now we go on to the general case. It is sufficient to show that \[\Delta_{0+}(v_0^rL_{0+}b_+^r) = (v_0^r\otimes v_0^r)\cdot (\sum_{p=0}^{\infty} (q^2;q^2)_p^{-1} \;v_0^pL_{0+}b_0^p\otimes v_0^pL_{0+}(-qb_+^*)^p)\cdot \Delta_+(b_+^r).\] Note then that \[v_0^rL_{0+} b_+^r= (q^2;q^2)_r^{1/2}\mathcal{G}_{0,r},\] so since $(\Delta_{0+}\otimes \iota)\mathcal{G}= \mathcal{G}_{[13]}\mathcal{G}_{[23]}$, we have \begin{eqnarray*} \Delta_{0+}(v_0^rL_{0+}b_+^r) &=&(q^2;q^2)_r^{1/2}\sum_t\mathcal{G}_{0,t} \otimes \mathcal{G}_{t,r} \\ &=& (q^2;q^2)_r^{1/2}\sum_t (q^2;q^2)_t^{-1/2} v_0^tL_{0+} b_+^t \otimes \mathcal{G}_{t,r},\end{eqnarray*} the sums converging $\sigma$-weakly. Now $\mathcal{G}_{t,r}$ is in $\mathbb{C}\lbrack \mathbf{X}\rbrack$ by the previous proposition. It is also not hard to see, using induction and the first commutation relation of Lemma \ref{LemCom1}, that the expression\[(v_0^r\otimes v_0^r)\cdot (\sum_{t=0}^{\infty} (q^2;q^2)_t^{-1}\, v_0^tL_{0+}b_+^t\otimes v_0^tL_{0+}(-qb_+^*)^t)\cdot \Delta_+(b_+^r)\] can be written in the form \[\sum_t (q^2;q^2)_t^{-1/2}\, v_0^tL_{0+}b_+^t \otimes f_{t,r},\] with $f_{t,r}\in \mathbb{C}\lbrack \mathbf{X}\rbrack$ and the sum $\sigma$-weakly converging again. We then have to prove that \[f_{t,r} = (q^2;q^2)_r^{1/2}\;  \mathcal{G}_{t,r}.\]

\noindent For this, it is enough to prove that \[\langle f_{t,r},x\rangle = (q^2;q^2)_r^{1/2}  \; \langle \mathcal{G}_{t,r},x\rangle \qquad \textrm{for all }x\in U_q(0,+).\] Now for each $y\in U_q(0,+)$ the pairing $\langle a_0^t\vartheta_{+0}^*b_+^t,y\rangle $ is only non-zero for finitely many $t$, so both \[\sum_t (q^2;q^2)_t^{-1/2}\, a_0^t\vartheta_{+0}^*b^t \otimes f_{t,r}\quad \textrm{and}\quad(q^2;q^2)_r^{1/2}\sum_t\mathcal{G}_{0,t} \otimes \mathcal{G}_{t,r}\] make sense as weakly convergent sums in the vector space dual $(U_q(0,+)\odot U_q(0,+))^{\circ}$. Since there exist elements $y_s \in U_q(0,+)$ such that \[\langle a_0^t\vartheta_{+0}^*b_+^t, y_s \rangle = \delta_{s,t},\] it is then sufficient to prove that the above expressions are equal as functionals on $U_q(0,+)\odot U_q(0,+)$. \\

\noindent But for $x,y\in U_q(0,+)$, we have  \begin{eqnarray*} (q^2;q^2)_r^{1/2}\, \langle \sum_t\mathcal{G}_{0,t} \otimes \mathcal{G}_{t,r}, x\otimes y\rangle &=& (q^2;q^2)_r^{1/2} \,\sum_t\langle \mathcal{G}_{0,t},x\rangle\,\langle \mathcal{G}_{t,r},y\rangle \\ &=& (q^2;q^2)_r^{1/2} \,\sum_t \Theta(x)_{0,t}\Theta(y)_{t,r} \\ &=&  (q^2;q^2)_r^{1/2} \,\Theta(xy)_{0,r}\\ &=& \langle a_0^r\vartheta_{+0}^*b_+^r,xy\rangle\end{eqnarray*}
on the one hand, while, since multiplication is separately continuous in each of its arguments in the weak topology,

\begin{eqnarray*} \langle \sum_t a_0^t\vartheta_{+0}^*b_+^t \otimes f_{t,r},x\otimes y\rangle &=& \langle (a_0^r\otimes a_0^r)\cdot (\sum_{t=0}^{\infty} (q^2;q^2)_t^{-1}\, a_0^t\vartheta_{+0}^*b_+^t\otimes a_0^t\vartheta_{+0}^*(-qb_+^*)^t)\cdot \Delta_+(b_+^r),x\otimes y\rangle \\ &=&   \langle \Delta_0(a_0^r)\Delta_{0+}(\vartheta_{+0}^*)\Delta_+(b_+^r),x\otimes y\rangle\\ &=& \langle \Delta_{0+}(a_0^r\vartheta_{+0}^*b_+^r),x\otimes y\rangle \\ &=& \langle a_0^r\vartheta_{+0}^*b_+^r,xy\rangle.\end{eqnarray*}

\noindent This concludes the proof.

\end{proof}

\noindent \emph{Remark:} The operator $\Delta_{0+}(L_{0+})$ also appears implicitly in \cite{Wor6}: the operator $X$ defined in formula (29) there is precisely $(|L_{0+}|^{-1/2}\otimes |L_{0+}|^{-1/2})\Delta_{0+}(L_{0+})^*$. In fact, another formula for this operator is given there, which could be used to give a more direct proof that this element satisfies a coassociativity condition.\\

\noindent We will relegate the closer study of the unitary $\mathcal{G}$ to a separate paper, since this belongs to the representation theory of $\mathbb{L}$ (in this paper we study rather (part of) the representation theory of $\mathbf{L}$). However, for completeness, we state here without proof the explicit form of the other entries of $\mathcal{G}$.\\

\noindent Denote by $p_n(x;a,0\mid q)$ the Wall polynomial of degree $n$ with parameter value $a$; so \[p_n(x;a,0\mid q) = \,_2\varphi_1(q^{-n},0;qa\mid q,qx),\] using the standard notation for $q$-hypergeometric functions (see \cite{Gas1}). For $t\leq s$, we then have \[\mathcal{G}_{t,s} = q^{t(t-s)}\left(\frac{(q^2;q^2)_s}{(q^2;q^2)_t}\right)^{1/2} (q^2;q^2)_{s-t}^{-1} \cdot v_0^{s+t} L_{0+} b_+^{s-t}\cdot  p_t(b_+^*b_+;q^{2(s-t)},0\mid q^2),\] while for $t\geq s$, \[\mathcal{G}_{t,s} = q^{s(s-t)}\left(\frac{(q^2;q^2)_t}{(q^2;q^2)_s}\right)^{1/2} (q^2;q^2)_{t-s}^{-1}\cdot v_0^{t+s} L_{0+} (-qb_+^*)^{t-s}  \cdot p_s(b_+^*b_+;q^{2(t-s)},0\mid q^2).\]

\noindent The proof can be carried out either by comparing the second factors in the expressions $\Delta_{0+}(\mathcal{G}_{0,s}) = \sum_t \mathcal{G}_{0,t}\otimes \mathcal{G}_{t,s}$ by means of Theorem \ref{TheoComu}, or by using Proposition \ref{PropVal}.\\

\noindent We remark that it is not so surprising to see the Wall polynomials appear here. Indeed, the projective $\mathbf{X}$-representation $\mathcal{G}$ can be seen as `the square root' of the coaction $\alpha$, which is in its form very similar to the $^*$-representation $(\theta\otimes \theta)\Delta$ of $C(SU_q(2))$ (using the notation from the proof of Proposition \ref{PropRow}). It is known that the Clebsch-Gordan coefficients of this representation are expressible precisely in terms of Wall polynomials (see \cite{Koo1}, Remark 4.2). For us however, these polynomials arise as the Clebsch-Gordan coefficients $\langle \mathcal{G}_{r,s}e_n^{+},e_m^{0+}\rangle$ of the $^*$-representation $\alpha$ of $C(S_{q0}^2)$. But both results essentially hinge on information concerning the spectrum of the operator $\Delta_+(b_+^*b_+)$.\\

\section{A Morita equivalence between the duals of $SU_q(2)$ and $\widetilde{E}_q(2)$}

\noindent In this section, we identify $\mathbb{H}$ with $\widetilde{E}_q(2)$.

\begin{Def} We define $n_0\, \eta \,\mathscr{L}(\mathbf{H})=\mathscr{L}^{\infty}(\mathbb{H})$ as the closure of \[\mathbb{C}\lbrack \mathbb{Z}\rbrack \odot \mathbb{C}\lbrack\mathbb{N}\rbrack\rightarrow \mathscr{H}_{0+}:e_n^{0+}\otimes e_k^{0+} \rightarrow q^n\;e_n^{0+}\otimes e_{k+1}^{0+}.\] Then $n_0$ is an unbounded normal operator.\end{Def}

\begin{Prop} The following formula holds: \[\Delta_{\mathscr{L}^{\infty}(\mathbb{H})}(v_0) = v_0\otimes v_0.\]\end{Prop}

\begin{proof} This is immediate from Theorem \ref{TheoComu} and Proposition \ref{PropDens}.\end{proof}

\begin{Prop} The following formula holds: \[\Delta_{\mathscr{L}^{\infty}(\mathbb{H})}(n_0) = n_0\otimes v_0 \dot{+} v_0^*\otimes n_0,\] with the latter being the closure of the sum of the two operators involved.\end{Prop}

\noindent \emph{Remark:} The computation which follows will be very similar to an argument appearing in \cite{Wor6}. We will also use the known properties of $n_0\otimes v_0 \dot{+} v_0^*\otimes n_0$, for example its normality (\cite{Wor4}).

\begin{proof} \noindent In the proof, we drop the sign index of basis vectors to lighten notation.\\

\noindent We first remark that it is sufficient to prove that \[(n_0\otimes v_0 \dot{+} v_0^*\otimes n_0)\Delta_{0+}(L_{0+}) = \Delta_{0+}(L_{0+})\Delta_+(b_+).\] Indeed, if this holds, then by the trivial commutation relation $n_0L_{0+} = L_{0+}b_+$, we get \[(n_0\otimes v_0 \dot{+} v_0^*\otimes n_0)\Delta_{0+}(L_{0+}) = \Delta_{0}(n_0)\Delta_{0+}(L_{0+}).\] Hence, using $v_0n_0v_0^*= qn_0$ and the previous proposition, we obtain \[\Delta_{0}(n_0 \chi_n) = (n_0\otimes v_0 \dot{+} v_0^*\otimes n_0)\Delta_{0}(\chi_n) \qquad \textrm{for all }n\in \mathbb{Z},\] where $\chi_n = \sum_{k=n}^{\infty} (e_{kk}\otimes 1)$.\\

\noindent But $|n_0|\chi_n\rightarrow |n_0|$ in the strong resolvent topology. Hence \[\Delta_{0}(n_0) \subseteq n_0\otimes v_0 \dot{+} v_0^*\otimes n_0,\] and since both operators are normal, equality follows.\\

\noindent We now prove the commutation relation \[(n_0\otimes v_0 \dot{+} v_0^*\otimes n_0)\Delta_{0+}(L_{0+}) = \Delta_{0+}(L_{0+})\Delta_+(b_+).\] Since $v_0n_0v_0^*= qn_0$, we have, for $p\in \mathbb{N}$, \[\begin{array}{l}(n_0\otimes v_0) \sum_{t=0}^p  (q^2;q^2)_t^{-1}v_0^t  L_{0+}b_+^t\otimes v_0^tL_{0+}(-qb_+^*)^t \\\\ \qquad \qquad =\; \sum_{t=0}^p (-1)^t (q^2;q^2)_t^{-1} v_0^t L_{0+}b_+^{t+1}\otimes v_0^{t+1}L_{0+}(b_+^*)^t.\end{array}\] In the same way, one computes that \[\begin{array}{l}(v_0^*\otimes n_0) \sum_{t=0}^p (q^2;q^2)_t^{-1}v_0^t L_{0+}b_+^t\otimes v_0^tL_{0+}(-qb_+^*)^t  \\ \\ \qquad \qquad =\;  \sum_{t=0}^p (-1)^t (q^2;q^2)_t^{-1} v_0^{t-1} L_{0+}b_+^{t}\otimes v_0^{t}L_{0+}(b_+^*)^tb_+.\end{array}\] Since the graph of $n_0\otimes v_0 \dot{+} v_0^*\otimes n_0$ is closed in the weak topology, it is then sufficient to show that \[\begin{array}{l} \sum_{t=0}^p (-1)^t (q^2;q^2)_t^{-1}v_0^t L_{0+}b_+^{t+1}\otimes v_0^{t+1}L_{0+}(b_+^*)^t  \\ \\ \qquad \qquad +\; \sum_{t=0}^p (q^2;q^2)_t^{-1}(-1)^t v_0^{t-1} L_{0+}b_+^{t}\otimes v_0^{t}L_{0+}(b_+^*)^{t}b_+ \end{array}\] converges to $\Delta_{0+}(L_{0+})\Delta(b_+)$ in the weak operator topology.\\

\noindent Write $x= q^{-1}b_+(e_{00}\otimes 1)$. Then $x$ is a bounded operator, and on $(\l^2(\mathbb{N})\otimes l^2(\mathbb{Z}))\otimes (l^2(\mathbb{N}_0)\otimes l^2(\mathbb{Z}))$, we have the above sum equals \[\begin{array}{l} \sum_{t=0}^p (-q)^t (q^2;q^2)_t^{-1} v_0^t L_{0+}b_+^{t+1}\otimes v_0^{t+1}L_{0+}(x^*)^t \\ \\ \qquad \qquad +\; \sum_{t=0}^p (-1)^t (q^2;q^2)_t^{-1} q^{t+1}v_0^{t-1} L_{0+}b_+^{t}\otimes v_0^{t}L_{0+}(x^*)^{t}x \end{array},\] which clearly converges in norm. Similarly, on $(l^2(\mathbb{N}_0)\otimes l^2(\mathbb{Z}))\otimes (l^2(\mathbb{N})\otimes l^2(\mathbb{Z}))$, we have that the sum equals \[\begin{array}{l}\sum_{t=0}^p (-1)^t (q^2;q^2)_t^{-1}q^{t+1}v_0^t L_{0+}x^{t+1}\otimes v_0^{t+1}L_{0+}(b_+^*)^t \\ \\ \qquad \qquad +\; \sum_{t=0}^p (-1)^t(q^2;q^2)_t^{-1} q^{t}v_0^{t-1} L_{0+}x^{t}\otimes v_0^{t}L_{0+}(b_+^*)^{t}b_+,\end{array}\] which again converges in norm. Hence, we only have to check the convergence (and the operator identity claim) on vectors of the form $(e_k\otimes \xi)\otimes (e_l\otimes \eta)$.\\

\noindent Now if $k=l=0$, one computes that the evaluation of the above sum on such a vector gives $(q^2;q^2)_{\infty}$ times the vector \[\begin{array}{l} e_1\otimes \xi \otimes e_0 \otimes S^{-1}\eta\\ \\ \qquad \qquad +\; \sum_{t=0}^{p-1} (-1)^t\, ((q^2;q^2)_t^{-1}-(q^2;q^2)_{t+1}^{-1})\;(e_{-t}\otimes S^{t+1}\xi\otimes e_{-(t+1)}\otimes S^{-t}\eta)\\ \\ \qquad \qquad \qquad \qquad \qquad \qquad \qquad +\; (-1)^p\; e_{-p}\otimes S^{p+1}\xi\otimes e_{-(p+1)} \otimes S^{-p}\eta.\end{array}\] On the other hand, evaluating \[(\sum_{t=0}^{p} (q^2;q^2)_t^{-1} v_0^t L_{0+}b_+^t \otimes v_0^tL_{0+}(-qb_+^*)^t)\cdot \Delta_+(b_+)\] on such a vector, we get $(q^2;q^2)_{\infty}$ times the sum of \[\begin{array}{l} e_1\otimes \xi \otimes e_0 \otimes S^{-1}\eta \\ \\ \qquad\qquad +\;\sum_{t=0}^{p-1} (-1)^{t+1} (q^2;q^2)_{t+1}^{-1} \,q^{2(t+1)}\; (e_{-t}\otimes S^{t+1}\xi \otimes e_{-(t+1)} \otimes S^{-t}\eta)\end{array}.\] The difference between the two evaluated vectors is \[(q^2;q^2)_{\infty}\,(-1)^p e_{-p} \otimes S^{p+1}\xi \otimes e_{-(p+1)} \otimes S^{-p}\eta.\] This indeed converges weakly to zero.\\

\noindent The case $kl\neq0$ is less subtle, since we have norm-convergence: we can then simply express both sides with respect to basis vectors, and compute that the resulting coefficients equal. We omit the details.\\

\end{proof}

\begin{Theorem} The locally compact quantum groups $\mathbb{H}$ and $\widetilde{E}_q(2)$ are isomorphic.\end{Theorem}

\begin{proof} This is immediate by the previous two propositions and the fact that $\mathscr{L}(\mathbf{H})$ is generated by $v_0$ and (the spectral projections of) $n_0$. (For the definition of $\widetilde{E}_q(2)$, see \cite{Wor2}, and \cite{Jac1} for a treatment in the setting of locally compact quantum groups.)
\end{proof}

\noindent \emph{Remarks:} \begin{enumerate}\item In fact, one can then also see this as a \emph{construction} of $\widetilde{E}_q(2)$ as a locally compact quantum group, albeit in a rather indirect way. This provides then a new, non-constructive proof of the existence of Haar weights on $\mathscr{L}^{\infty}(\widetilde{E}_q(2))$ (which was already proven in two quite different ways in \cite{Baa2} and \cite{Pal1}).
\item The action of $SU_q(2)$ on the standard Podle\'{s} sphere of course restricts to $SO_q(3)$, and then the above restriction will give a Morita equivalence between the duals of $SO_q(3)$ and $E_q(2)$. While it is valid to say that this is the more natural setting, we prefer to work with the 2-folded coverings since the computations are somewhat easier.
\end{enumerate}

\noindent We can now easily prove the statement we made in the abstract.

\begin{Theorem}\label{TheoCoc} There exists a unitary 2-cocycle $\Omega\in \mathscr{L}^{\infty}(SU_q(2))\bar{\otimes}\mathscr{L}^{\infty}(SU_q(2))$ such that \[(\mathscr{L}^{\infty}(SU_q(2)),\Delta_{\Omega}) \cong \mathscr{L}^{\infty}(\widetilde{E}_q(2)).\]
\end{Theorem}

\begin{proof} Let $\tilde{x}$ be a unitary $l^2(\mathbb{N})\rightarrow l^2(\mathbb{Z})$, and put \[x=\tilde{x}\otimes 1: \mathscr{H}_+\rightarrow \mathscr{H}_{0+}.\] Then $x\in \mathscr{L}(\mathbf{X})$, and we put \[\Omega := (x^*\otimes x^*)\Delta_{0+}(x).\] It is then straightforward to check that $\Omega$ is a unitary 2-cocycle. The map \[\mathscr{L}^{\infty}(\widetilde{E}_q(2)) \rightarrow \mathscr{L}^{\infty}(SU_q(2)): y\rightarrow x^*yx\] provides the isomorphism stated in the theorem.

\end{proof}

\noindent\emph{Remarks:} \begin{enumerate} \item This shows that the regularity of the standard left multiplicative unitary of a von Neumann algebraic quantum group (\cite{Baa1}) is not preserved by cocycle twisting, since it is well-known that $\widetilde{E}_q(2)$ is not regular. It is an interesting question to see if there are necessary and sufficient conditions for regularity to be preserved.
\item In \cite{Wor6}, Woronowicz proves that, with $u:\mathscr{H}_+\rightarrow \mathscr{H}_{0+}$ the canonical isometry as before, the formula $\Delta_{0}(uxu^*) = Z^* \Delta_+(x)Z$ holds for some \emph{isometry} $Z \in \mathscr{L}(\overline{\mathbf{X}})\bar{\otimes}\mathscr{L}(\overline{\mathbf{X}})$, which (in his approach) is easily seen to satisfy the 2-cocycle relation. In fact, it is also easily seen that, in our notation, $Z= \Delta_{+0}(u^*)$. It is not clear to us if Woronowicz' methods could be refined as to provide, in a straightforward fashion, a \emph{unitary} 2-cocycle as in the previous theorem, but we should admit that we have not considered this possibility in detail.
\end{enumerate}

\vspace{0.1cm}

\noindent \emph{Acknowledgements}: I would like to thank Pawe{\l} Kasprzak and Wojciech Szyma\'{n}ski for their invitation to Syddansk Universitet, Odense, Denmark, where part of this work was done. I would also like to thank Erik Koelink and Julien Bichon for providing some extra references.

\end{document}